\documentclass[10pt]{amsart}

\usepackage{xypic,amsfonts,amsmath,amssymb}
\usepackage{pstricks,pst-node,pst-plot}
\usepackage{multicol}

\CompileMatrices

\xyoption{all}
\input xypic
\newtheorem{theorem}{Theorem}[section]
\newtheorem{proposition}[theorem]{Proposition}

\newtheorem{definition}[theorem]{Definition}
\newtheorem{lemma}[theorem]{Lemma}
\newtheorem{remark}[theorem]{Remark}

\newcommand{\pf}{\noindent \emph{\mbox{Proof:\ \ \ }}}
\newcommand{\cvd}{\hfill $\Box$}

\newcommand{\ie}{{\it i.e.\/}\ }

\newcommand{\ra}{\rightarrow}

\newcommand{\arup}[1]{\stackrel{#1}{\rightarrow}}
\newcommand{\arr}[1]{\stackrel{#1}{\longrightarrow}}

\newcommand{\iso}{\stackrel{_\sim}{\rightarrow}}
\newcommand{\epi}[1]{\stackrel{#1}{\twoheadrightarrow}}
\newcommand{\mono}[1]{\stackrel{#1}{\hookrightarrow}}

\newcommand{\cd}{{\mathcal D}}
\newcommand{\ce}{{\mathcal E}}
\newcommand{\cf}{{\mathcal F}}
\newcommand{\cg}{{\mathcal G}}
\newcommand{\ch}{{\mathcal H}}

\newcommand{\ck}{{\mathcal K}}
\newcommand{\cl}{{\mathcal L}}
\newcommand{\cm}{{\mathcal M}}

\newcommand{\co}{{\mathcal O}}

\newcommand{\cq}{{\mathcal Q}}

\newcommand{\ct}{{\mathcal T}}
\newcommand{\cu}{{\mathcal U}}
\newcommand{\cv}{{\mathcal V}}

\newcommand{\Z}{{\mathbf{Z}}}

\newcommand{\R}{{\mathbf{R}}}
\newcommand{\N}{{\mathbf{N}}}

\newcommand{\ts}{\otimes}
\newcommand{\Ts}{\bigotimes}

\newcommand{\Hom}{{\mbox{\rm{Hom}}}}
\newcommand{\im}{\mbox{Im}}
\newcommand{\coker}{\mbox{\mbox{\rm coker}\ }}

\newcommand{\fhom}{\mathop{\mathcal{H}om}}
\newcommand{\ext}{{\mbox{\rm{Ext}}^1}}
\newcommand{\rk}{\mathop{\rm rk}\nolimits}
\newcommand{\grass}{\mathop{\rm Grass}\nolimits}
\newcommand{\flag}{\mathop{\rm Flag}\nolimits}
\newcommand{\Flag}{\mathop{\mathcal{F}\rm {lag}}\nolimits}
\newcommand{\quot}{\mathop{\rm Quot}\nolimits}
\newcommand{\qpquot}{\mathop{{\rm Quot}^{qpar}_{(n''_i),d''}}\nolimits}
\newcommand{\pquot}{\mathop{{\rm Quot}^{par}_{r'',s_*}}\nolimits}
\newcommand{\qpquotI}{\mathop{{\rm Quot}^{qpar}_{((n''_i(p))_{p\in I}),d''}}\nolimits}
\newcommand{\qpquoti}{\mathop{{\rm Quot}^{qpar}_{(n''_i(p)),d''}}\nolimits}
\newcommand{\QI}{\mathop{{\rm Q}_{((n'_i(p))_{p\in I}),d''}}\nolimits}
\newcommand{\stack}{\mathop{{\mathcal QPar}_{X,I,\underline n}}\nolimits}

\newdir{ >}{{}*!/-8pt/\dir{>}}

\title{\bf Theta functions on the moduli space of parabolic bundles}

\date\today

\begin{document}

\author{Francesca Gavioli}
\address{Laboratoire de Math\'ematiques Jean Leray CNRS UMR 6629\\ 
Facult\'e des Sciences, Universit\'e de Nantes\\ 
2 rue de la Houssini\`ere, BP 92208\\ 
F 44322 Nantes Cedex 03, France}
\email{francesca.gavioli@math.univ-nantes.fr}

\maketitle

\markboth{\author}{Parabolic theta functions}

\section{\bf Introduction}

\bigskip

Let $X$ be a smooth, connected projective curve of genus $g \ge 2$
over the field of complex numbers and $I$ a finite subset of points of
$X$. Let $\cm^{par}$ denote the moduli space of semistable parabolic vector
bundles of rank $r$, trivial determinant and fixed parabolic
structure at $I$. There is a natural ample line bundle $\cl^{par}$ on
$\cm^{par}$, which is the analogue of the determinant bundle $\cd$ on the
moduli space $\cm$, of vector bundles on $X$ of fixed rank and
determinant (\cite{narasimhan-ramadas} theorem 1, for rank $2$, 
\cite{pauly 1} theorem 3.3, for any rank). In this paper we determine
an integer $\ell_0$ such that, if $\ell \ge \ell_0$, then 
${\cl^{par}}^{\ts \ell}$ is globally generated.\\

The analogue problem in the classical case has been studied by
Faltings, Le Potier and Popa. For vector bundles,
there are natural global sections (of each power $h$) of the
determinant bundle, that are called theta functions (of order
$h$). Faltings has shown that such sections do generate $\cd^{\ts h}$, for
$h \gg 0$ \cite{faltings}, and an effective bound on $h$ has been given
by Le Potier \cite{jlp}. Recently Popa has produced a considerably better
bound, in the sense that it does not depend on the genus $g$ of the curve 
\cite{popa}. \\

The parabolic case in rank $2$ has been studied by Pauly \cite{pauly 2}. 
He produces
sections of the parabolic determinant $\cl^{par}$ on the moduli space of 
semistable parabolic bundles of rank $2$ and trivial determinant. 
They generalize the sections of type theta of the determinant 
line bundle. 
Moreover, under the assumption that the parabolic subset $I$ 
has small and even cardinality, he proves that these sections generate the 
line bundle $\cl^{par}$.\\

Our main result can be stated as follows. \\

\begin{theorem}\label{gg} Let $\ell$ be an integer such that
$$
\ell \ge \left[\frac{r^2}{4}\right],
$$ 
and suppose it is $I \neq \emptyset$. 
Then the linear system $|{\cl ^{par}}^{\ts\ell}|$ 
is base point free.\\
\end{theorem}

\medskip

We are actually going to prove that, for $\ell$ given by this bound,
there exist global sections, the parabolic analogues of theta
functions, generating ${\cl ^{par}}^{\ts\ell}$. These sections are obtained
generalizing Pauly's method \cite{pauly 2} and will be
called parabolic theta functions. 
They are associated with parabolic bundles whose rank, degree and parabolic
invariants depend on the invariants of the bundles parametrized 
by $\cm^{par}$ and on the order $\ell$. 
Let $\cm '_\ell$ denote the moduli space of semistable parabolic
bundles with which we associate parabolic theta functions of order $\ell$. \\

The idea of the proof is to show that, under this assumption on
$\ell$, for each point $x$ of the
moduli space $\cm^{par}$, the dimension of the subscheme of points of
$\cm'_\ell$, whose associated parabolic theta function vanishes at $x$, is
strictly smaller than $\dim(\cm'_\ell)$. \\

Our method of proof is inspired by Popa's beautiful ideas \cite{popa}.
An essential step in his proof is the estimate 
of the dimension of Grothendieck's Quot scheme (see also 
\cite{popa-roth}). This allows him to estimate the dimension of
the family of bundles, that are images of a morphism from a fixed 
vector bundle $E$. \\ 
In order to treat the parabolic case, given a point
$F_*$ of $\cm'_\ell$, we first show how to identify the
zeroes of the associated section with the points $E_*$ of
$\cm^{par}$, admitting a nonzero parabolic morphism to $F_*$. 
Particular care has then to be taken in order to understand the family of 
such morphisms for which we construct a scheme that can be seen as
a parabolic analogue of Grothendieck's Quot scheme. 
An estimate of the dimension of this scheme is given by a formula 
(see theorem \ref{cong}), which extends
the result of Popa and Roth \cite{popa-roth} to the parabolic case.
The computation can then be worked out, by applying
Lange's results on families of extensions to the parabolic context. \\

\bigbreak

Let $\stack$ be the (algebraic) stack of quasi-parabolic vector bundles
with trivial determinant and fixed quasi-parabolic structure. 
It is well known that the choice of a system $\underline \alpha$ of Seshadri 
parabolic weights defines a notion of $\underline \alpha$-semistability 
for such bundles. 
Actually by the classification of line bundles on $\stack$ of Pauly 
\cite{pauly 1} and Laszlo and Sorger \cite{laszlo-sorger} 
the choice of $\underline \alpha$ also defines a line bundle 
$L_{\underline \alpha}$ on this stack and 
all ``ample'' line bundles arise this way.
The above theorem may be applied to show that for $\ell$ sufficiently
large, the base locus of the linear system 
$|L_{\underline \alpha}^{\ts \ell}|$ on $\stack$ is isomorphic to the closed 
substack of $\underline \alpha$-unstable parabolic bundles. 
%In this sense, it turns out that
%Seshadri's notion of semistability involving parabolic weights is the
%``only'' reasonable notion of semistability for quasi-parabolic bundles.\\

\bigskip

{\small {\bf Acknowledgements.} I am deeply grateful to  
Christoph Sorger for the interest he has taken in this work and
to Christian Pauly for several interesting and helpful discussions.}\\

\bigskip

\section{\bf Parabolic bundles}

\bigskip

A {\it quasi-parabolic bundle} $(E, (f_p)_{p \in I})$ on $X$ with
quasi-parabolic structure at $I$ is a vector bundle $E$ on $X$ and
flags $f_p$ of the fibre of $E$ over $p$, for $p\in I$:
$$
E_{p} = E_{p,1} \supset E_{p,2} \supset \dots \supset E_{p,l_{p}}
\supset E_{p,l_{p} +1}=0.
$$
The positive integers $n_i(p)=\rk(E_{p,i}/E_{p,i+1})$ are the {\it
  multiplicities} of $(E,(f_p)_{p \in I})$ at $p$ and $l_p$ is the
  {\it length of the flag $f_p$}. Let $r_i(p)$ denote
  $\sum_{j=1}^in_i(p)$.\\

This is equivalent to considering filtered locally free sheaves
$$
E= E _{(p,1)} \supset E _{(p,2)} \supset \dots \supset E _{(p,l_p)}
\supset E _{(p,l_p+1)} = E(-p),
$$
where $E_{(p,i)} =\ker(E\twoheadrightarrow E_p \twoheadrightarrow
E_p/E_{p,i})$. This filtration will again be denoted by $f_p$. Here
the multiplicities are defined as the integers $n_i(p)=\deg(E
_{(p,i)}/E _{(p,i+1)})$.\\  Let $\flag_{n_1, \dots , n_l}(E_p)$ be the
flag variety of $E_p$ of type $(n_1, \dots , n_l)$. It is an
irreducible projective variety of dimension $d_{n_1, \dots ,
n_l}=\displaystyle\sum_{i>j}n_in_j$. \\

Let $E \epi q G$ be a quotient bundle of $E$. Then a quasi-parabolic
structure on $E$ induces a quasi-parabolic structure on $G$: let
$h_{p,i}$ be the injection $E_{(p,i)} \hookrightarrow E$ and denote
$G_{(p,i)}=\im(qh_{p,i})$. Then the quotient morphism induces a
filtration at each parabolic point
$$
G= G_{(p,1)} \supseteq G_{(p,2)} \supseteq \dots \supseteq G_{(p,l_p)}
\supseteq G_{(p,l_p+1)}=G(-p).
$$
By considering the distinct locally free sheaves of each filtration,
this defines a  quasi-parabolic structure on $G$. It is induced by the
one on $E$ in the sense that the morphism $q$ is naturally compatible
with the filtrations. Dually, there is a natural induced
quasi-parabolic structure on a subbundle $H \mono j E$: if
$\pi=\coker(j)$, then it is obtained by letting $H_{(p,i)}= \ker(\pi
h_{p,i})$. In other words, $H_{(p,i)}=H \cap E_{(p,i)}$.\\

Let $(V'',({f_p}_{V''})_{p \in I})$ (respectively, $(V',(f_p^{V'})_{p
\in I})$) denote the quasi-parabolic structure induced by $(E,(f_p)_{p
\in I})$ on a quotient bundle $V''$ (respectively, a subbundle $V'$).\\

\bigskip

A {\it parabolic bundle} $E_*$ on $X$ is a quasi-parabolic bundle
with, for all $p \in I$, a sequence of real numbers
$$
0 \le \alpha_1(p) < \alpha_2(p) < \dots < \alpha_{l_p}(p) <1,
$$ 
attached to the flag at $p$. These numbers are called {\it parabolic
weights}. It is convenient to introduce Simpson's equivalent definition
\cite{simpson} of a parabolic bundle as a filtered vector
bundle. In the notations of \cite{maruyama-yokogawa}, \cite{yokogawa},
a parabolic bundle is:

- for all $\alpha \in \R$, a locally free sheaf $E_\alpha$ on $X$ and
  an isomorphism 
$$
j_\alpha :E_\alpha (-\displaystyle\sum_{p \in I} p)
  \iso E_{\alpha +1},
$$

- for all $\alpha,\,\beta \in \R$, such that $\alpha \ge \beta$, an
  injective morphism $i_{E_*}^{\alpha, \beta} : E_\alpha
  \hookrightarrow E_\beta$, such that the diagram
$$\xymatrix{ E_{\alpha +1} \ar@{^{(}->}[r]^{\ \ i_{E_*}^{\alpha,
\alpha +1}}  & E_\alpha \\ E_\alpha(-\sum_{p \in I} p)
\ar[u]^{j_\alpha}\ar@{^{(}->}[r] & E_\alpha \ar[u]^{id}\\ }$$ commutes,

- a sequence of real numbers $0 \le \alpha_1 < \alpha_2 < \dots <
  \alpha_L <1,$ such that $i_{E_*}^{\alpha, \alpha_i}$ is an
  isomorphism $E_\alpha \cong E_{\alpha_i}$, for all $\alpha \in \
  ]\alpha_{i-1},\alpha_i]$. \\

As a convention for a parabolic bundle $E_*$ the sheaf $E=E_0$ is called 
the underlying vector bundle and for $\alpha \in \R^+$ the
morphisms of the parabolic structure will be denoted by 
$i_E^\alpha \colon= i_{E_*}^{\alpha, 0}$ and
$\pi_E^\alpha \colon= \coker(i_E^\alpha)$. \\

Let $E_*$ and $F_*$ be parabolic bundles on $X$, with parabolic
structure at $I$. A morphism $\varphi : E \ra F$ is parabolic if, for
all $\alpha\in \R^+$, the composition $\pi_F^\alpha \varphi
i_E^\alpha$ is the zero morphism. This produces a morphism
$$\xymatrix{ 0 \ar[r] & E_\alpha \ar[r]^{i_E^\alpha} \ar@{-->}[d] & E
\ar[r] \ar[d]^{\varphi} & E/E_\alpha \ar[r] & 0\\ 0 \ar[r] & F_\alpha
\ar[r] & F \ar[r]^{\pi_F^\alpha} & F/F_\alpha \ar[r] & 0\\ }$$ that
will be denoted by $\varphi_\alpha : E_\alpha \ra F_\alpha$. The notation
$\varphi_* : E_* \ra F_*$, means that $\varphi$ is parabolic.\\

Consider the sheaf defined by
$$
\fhom(E_*,F_*)(\cu)=\{\varphi_*^\cu : {E_*}_{|\cu} \ra {F_*}_{|\cu}\}.
$$
By definition of parabolic morphism, it
is a subsheaf of $\fhom(E,F)$ and for all open
subset $\cu \subset X$, such that $I \cap \cu=\emptyset$, it  actually
is $\fhom(E_*,F_*)(\cu)=\fhom(E,F)(\cu)$. Thus the quotient sheaf is
a torsion sheaf with support at $I$ that can be described in terms
of the parabolic structures of $E_*$ and $F_*$; this is done in  
\cite{boden-hu}, lemma 2.4: suppose for simplicity that $I=\{p\}$ and 
let $(\alpha_1,\dots,\alpha_l)$ be the weights,
$n_i^{E_*}=\deg(E_{\alpha_i}/E_{\alpha_{i+1}})$ the multiplicities of
$E_*$, $(\beta_1, \dots, \beta_h)$ the weights,
$n_j^{F_*}=\deg(F_{\beta_j}/F_{\beta_{j+1}})$ the multiplicities of
$F_*$. Then there is a short exact sequence
\begin{eqnarray}\label{parhom-hom-tau}\xymatrix{
0 \ar[r] & \fhom(E_*,F_*) \ar[r] & \fhom(E,F) \ar[r] &
\tau_{E_*,F_*}\ar[r] & 0,\\ }\end{eqnarray} where $\tau_{E_*,F_*}$ is
a torsion sheaf supported at $p$ of degree
$$
h^0(X,\tau_{E_*,F_*})=\sum_{{i,j}\atop{\alpha_i >
\beta_j}}n_i^{E_*}n_j^{F_*}.
$$
This is a consequence of the fact that a morphism $\varphi : E \ra F$
is parabolic, if and only if the linear map over the parabolic point
$$
\varphi _p =(\varphi _{i,j}) : \bigoplus _i
  E_{\alpha_i}/E_{\alpha_{i+1}} \arr {} \bigoplus _j
  F_{\beta_j}/F_{\beta_{j+1}}
$$
is such that $\varphi _{i,j}=0$, for all $\alpha_i>\beta_j$. Hence the
fiber of $\tau_{E_*,F_*}$ at $p$ is isomorphic to
$$
\bigoplus _{\alpha_i>\beta_j}(E_{\alpha_i}/E_{\alpha_{i+1}})^\vee \ts
F_{\beta_j}/F_{\beta_{j+1}}.
$$
For a general parabolic subset $I$, the degree of the torsion sheaf
$\tau_{E_*,F_*}$ can be computed as
$$h^0(X,\tau_{E_*,F_*})=\sum _{p  \in I}\sum_{{i,j}\atop{\alpha_i(p) >
\beta_j(p)}}n_i^{E_*}(p)n_j^{F_*}(p) .
$$
Let $\chi(E_*,F_*)$ denote $\chi(\fhom(E_*,F_*))$. By Riemann-Roch
formula and the exact sequence (\ref{parhom-hom-tau}), this Euler
characteristic can be computed as
$$
\chi(E_*,F_*)=\rk(E) \deg(F) -\rk(F)\deg(E) +\rk(E)\rk(F)(1-g)-
h^0(\tau_{E_*,F_*}).
$$
The group of global sections $H^0(\fhom(E_*,F_*))$ is the group of
parabolic morphisms, $\Hom(E_*,F_*)$. The first cohomology group
$H^1(\fhom(E_*,F_*))$ is, by \cite{yokogawa} lemma 1.4, isomorphic
to the group of isomorphism classes of parabolic extensions of $F_*$
by $E_*$ and is denoted by $\ext(E_*,F_*)$. By definition, a parabolic
extension is a short exact sequence
$$\xymatrix{ 0 \ar[r] & E '_* \ar[r]^{i_*} & E_* \ar[r]^{p_*} & E ''_*
\ar[r] & 0, }$$ two parabolic extensions being isomorphic, if there is
a parabolic isomorphism of extensions.\\

Recall the definitions of the parabolic invariants of $E_*$. Let $d=
\deg(E)$ and $r=\rk(E)$; the {\it parabolic degree} of $E_*$ is
defined as the real number
$$
\deg(E_*)=\deg(E) +\sum_{p\in I} \sum _{i=1}^{l_p} n_i(p) \alpha
_i(p) 
$$
and can be computed as the integral
$$
\int _{0}^1\deg(E_\alpha)d\alpha + r|I|=\int _{-1}^0\deg(E_\alpha)d\alpha.
$$
The {\it parabolic Hilbert polynomial} is 
$$
{\rm P}(E(m)_*)=\deg(E(m)_*)+r(1-g)=\deg(E_*)+r(m+1-g)
$$ 
and the {\it parabolic slope} is $\mu(E_*)=\displaystyle\frac{\deg(E_*)}{r}$.\\

Let $E_*$ be a parabolic bundle and $E \epi {} G$ a quotient vector
bundle. Consider the parabolic structure obtained by the induced
quasi-parabolic structure, weighted as the one of $E_*$. This
parabolic structure is said to be the one induced on $G$ by $E_*$ and
will be denoted by $G_*$. Dually, all subbundle $H\mono {} E$ has an
induced parabolic structure, that will be denoted by $H_*$. Recall the
notations of the induced quasi-parabolic structure. For all $p \in I$
and $i= 1, \dots , l_p$, consider the integers $n''_i(p)=
\deg(G_{(i)}/G_{(i+1)})$. They verify $0 \le n''_i(p) \le
n_i(p)$ and $n''_1(p)+ \dots + n''_{l_p}(p)=\rk(G)$. Then we can
easily check the equality
$$
\deg (G_*)= \deg(G) + \sum_{p \in I} \sum_{i=1}^{l_p} n''_i(p) \alpha
_i(p). \\
$$

\begin{remark} \rm For all parabolic structure $G_{*'}$ (respectively, $H_{*'}$), such that $E_* \epi {} G_{*'}$ (respectively, $H_{*'} \mono {} E_*$) is parabolic, it is $\deg (G_*) \le \deg (G_{*'})$ (respectively, $\deg
  (H_{*}) \ge \deg (H_{*'})$).\\
\end{remark}

\begin{definition} A parabolic bundle $E_*$ is semistable if, for all quotient bundle $E \epi {} G$, the inequality $\mu(E_*)\le \mu(G_*)$ holds. A semistable bundle is stable if the inequality is strict, whenever $G$ is a nontrivial quotient of $E$.\\
\end{definition}

Suppose, to simplify the formulation of the following basic facts,
that $I=\{p\}$.\\

It is easily seen, that the (semi)stability of a parabolic bundle, as
a function of the weights, just depends on their differences. More
precisely, let $E_{*'}$ be the bundle with same quasi-parabolic
structure as $E_*$ and weights
$$
0 \le 0 < \delta_1 < \dots < \delta_1 + \dots + \delta_{l-1} < 1,
$$
where $\delta_i=\alpha_{i+1}-\alpha_i$. Then $E_*$ is (semi)stable if
and only if $E_{*'}$ is (semi)stable. Thus we can (and will) assume in
the following that the smallest weight at each parabolic point is
zero. With the notations of \cite{boden-hu} this means that we represent
the weights in the face $\partial_0 W$ of $W$. \\ 
It will be useful to
remark that this assumption allows to write the parabolic degree as
$$
\deg(E_*)=\deg(E)+\sum_{p\in I}\sum_{i>j}n_i(p)\delta_j(p).\\
$$

As it is shown in \cite{seshadri-mehta}, {\S} 2 the (semi)stability
condition actually depends on rational weights, \ie there is a
rational system of weights   $(\alpha'_1, \dots, \alpha'_l)$, such
that  $E_*$ is (semi)stable if and only if it is (semi)stable with
respect to the weights $(\alpha'_j)$. For this reason, in what follows
we consider rational weights.\\

With the notations of \cite{boden-hu} for the variation of the
(semi)stability condition, for a fixed quasi-parabolic structure
there exists an open subset of the space of weights
$W$ such that, for any system of
weights in this open subset, the condition of semistability is
equivalent to the condition of stability. 
Actually, this open subset is the complement of a union of hyperplanes 
\cite{seshadri-mehta}, that we will call {\it Seshadri walls}.
Writing the weights $(\alpha_i)$ as rational numbers $\alpha_i=
\frac{a_i}{k}$, these hyperplanes are given by the equations 
$h=(r, k\deg(E_*))$, for some integer $h \ge 2$.\\

\bigskip

\section{\bf {The schemes of quasi-parabolic and parabolic quotients}}

\bigskip

Let $(E,f)$ be a rank $r$ vector bundle endowed with a quasi-parabolic 
structure at $p$ of
multiplicities $(n_1, \dots, n_l)$. Consider the set of quotient
bundles $G$ of $E$, whose induced quasi-parabolic structure is of
fixed type $(n''_1, \dots, n''_l)$, that is if $(G,f_G)$ is the
induced structure, then $f_G \in \flag_{n''_1, \dots, n''_l}(G_p)$. In
fact this set can be equipped with a natural algebraic structure: in
the first part of this section we construct a subscheme of Grothendieck's 
scheme of quotients, parametrizing quotient bundles of fixed induced
multiplicities.\\
The same construction applies to the case of quotients of a parabolic
bundle $E_*$. This produces a scheme parametrizing quotient bundles of
fixed induced parabolic type.\\

Let $\quot_{r'',d''}(E)$ be the scheme of quotients of $E$ of rank
$r''$ and degree $d''$ and let $\quot^o_{r'',d''}(E)$ denote the open
subscheme of quotients $E \epi q G$, such that $G$ is a locally free
sheaf. Denote by $\pi_X$, $\pi_Q$ the projections of $X\times
\quot^o_{r'',d''}(E)$ on $X$ and $\quot^o_{r'',d''}(E)$ respectively
and $\pi_X^*(E) \epi \pi \cg$ the universal quotient. The morphism
$h_i : E_{(i)}\ra E$ produces an injective morphism $\pi_X^* h_i$ of 
locally free sheaves on $X\times \quot^o_{r'',d''}(E)$ and the image of the
composition $\pi \pi_X^*(h_i)$ is a subsheaf $\cg_{(i)}$ of the
universal family of quotients $\cg$:
$$
\xymatrix{  
\dots \ar[r] & \pi_X^* E_{(3)} \ar@{^{(}->}[r]^{\ \
\pi_X^*(h_3)} \ar@{>>}[d]^{\pi _3} & \pi_X^* E_{(2)}
\ar@{^{(}->}[r]^{\ \ \pi_X^*(h_2)} \ar@{>>}[d]^{\pi _2} & \pi_X^*E
\ar@{>>}[d]^{\pi}\\ \dots \ar[r] & \cg_{(3)} \ar@{^{(}->}[r]
&\cg_{(2)} \ar@{^{(}->}[r] & \cg \\ 
}.\\
$$ 
There is a flattening
stratification  of $\cg_{(2)}$ on $X\times \quot^o_{r'',d''}(E)$. Thus
we can write the scheme of locally free quotients as a disjoint union
$$
\quot^o_{r'',d''}(E)= \coprod_{\nu_1=0}^{r''} Q_{\nu_1}
$$
where, set theoretically, the scheme $Q_{\nu_1}$ consists of those 
quotients that can
be seen as points of $\quot_{r'',d''-\nu_1}(E_{(2)})$
via the induced surjective map:
$$
Q_{\nu_1}= \{[E \epi q G] \in \quot^o_{r'',d''}(E) | \deg
(G_{(2)})=d''-\nu_1\}. 
$$
Restrict the filtration of the universal family to the stratum
$Q_{\nu_1}$ and consider the sheaf $\cg_{(3)}$ as an $\co _{X \times
Q_{\nu_1}}$-module. Then there is a flattening stratification of
$Q_{\nu_1}$ with respect to the family $\cg_{(3)}$, hence we can write:
$$
Q_{\nu_1}=\coprod_{\nu_2=0}^{r''}Q_{\nu_1\nu_2}.
$$
Thus, taking flattening stratifications of each stratum, we end up
with a stratification of $Q_{\nu_1\dots\nu_{i-1}}$, with respect to
$\cg_{({i+1})}$:
$$
Q_{\nu_1\dots\nu_{i-1}}=\coprod_{\nu_i=0}^{r''}Q_{\nu_1\dots
\nu_{i-1}\nu_i}.\\
$$

\begin{remark} \rm Let $\nu_1, \dots , \nu_i$ be positive integers such that 
$\nu_1+ \dots + \nu_i>r''$. Then the stratum $Q_{\nu_1\dots\nu_i}$ is empty.\\
\end{remark}

This is straightforward, since there is a natural isomorphism
$\cg_{({l+1})}\iso \cg(-p)=\cg\otimes \pi_X^* \co _X(-p)$.\\ 

Thus, the induced quasi-parabolic type on quotients gives a
stratification of Grothendieck's scheme of locally free quotients
$$
\quot^o_{r'',d''}(E)= \coprod_{\nu_1+\dots + \nu_{l-1}=0}^{r''}
  Q_{\nu_1\dots\nu_{l-1}}.\\
$$

\begin{definition} Let $(n''_1, \dots, n''_l)$ be integers, such that 
$0\le n''_i \le n_i$ and $n''_1 + \dots + n''_l = r''$. We define the
scheme of quasi-parabolic quotients of $(E,f)$ of type $(n''_1, \dots,
n''_l)$ as the stratum $Q_{n''_1\dots n''_{l-1}}$ and will  denote it
by $\qpquot(E,f)$.\\
\end{definition}

Let $I\subset X$ be a finite subset of points and $(E,(f_p)_{p\in I})$ a
quasi-parabolic structure on $E$. Then quotient bundles of $E$ of rank
$r''$, degree $d''$ and fixed induced quasi-parabolic type are
parametrized by a locally closed subscheme of Grothendieck's scheme of
quotients. Let $((n''_i(p))_{p\in I})$ be integers such that, for all
$p \in I$, it is $0\le n''_i(p) \le n_i(p)$, for all $i=1, \dots ,
l_p$ and $\sum_i n''_i(p)=r''$. We define the {\it scheme of
  quasi-parabolic quotients} of $(E,(f_p)_{p\in I})$ of type $((n''_i(p))_{p\in
  I})$ as the intersection 
$$
\qpquotI(E,(f_p)_{p\in I})=\bigcap_{p\in I}\qpquoti(E,f_p).\\
$$

\bigskip

Let $E_*$ be a parabolic bundle. The construction of the scheme of quotients of $E_*$ of fixed induced parabolic structure is completely analogue to the construction of the scheme of quasi-parabolic quotients. It is enough to consider flattening stratifications of the families $\cg_{\alpha_i}$:
$$
\xymatrix{  
\dots \ar[r] & \pi_X^* E_{\alpha_2} \ar@{^{(}->}[r]^{\ \
\pi_X^*(i_{E_*}^{\alpha_2, \alpha_1})} \ar@{>>}[d]^{\pi _{\alpha_2}} & \pi_X^* E_{\alpha_1}
\ar@{^{(}->}[r]^{\ \ \pi_X^*(i_{E_*}^{\alpha_1})} \ar@{>>}[d]^{\pi _{\alpha_1}} & \pi_X^*E
\ar@{>>}[d]^{\pi}\\ \dots \ar[r] & \cg_{\alpha_2} \ar@{^{(}->}[r]
&\cg_{\alpha_1} \ar@{^{(}->}[r] & \cg \\ }.
$$
The strata are here determined by the Hilbert polynomials at the weights $\alpha_i$. This fixes the parabolic Hilbert polynomial of the quotients, since this
is determined by the Hilbert polynomials at each weight. \\

There are more useful multiplicities than the quasi-parabolic ones, 
that we introduce here, since they better fit to the parabolic filtration. 
For a point of a stratum, represented by a quotient bundle $G_*$, let $n''_{\alpha_1}$ be the positive integer defined by $\deg(G_{\alpha_1})= d'' - n''_{\alpha_1}$ and in general $n''_{\alpha_i}$ the integer such that $\deg(G_{\alpha_i})=\deg(G_{\alpha_{i-1}})-n''_{\alpha_i}$. 
In
analogy with the quasi-parabolic multiplicities, we will denote
$r''_{\alpha_i}=\sum_{j=1}^i n''_{\alpha_j}$.\\ 

By letting $\alpha_0=0$
the parabolic degree 
$d''_*$ can be viewed as a polynomial in the weights, 
in fact it can be computed as 
$$
d''_*= d'' + \sum_{i=1}^{L+1}n''_{\alpha_i}\alpha_{i-1}.\\
$$

The parabolic structure induced on a quotient bundle $G$ of 
$E_*$ is determined 
by the decreasing step function, that we denote by $s_*$, associated with 
the collection of 
the degrees at each weight, that is the {\it parabolic degree function} 
$$
\alpha \mapsto s_{\alpha}=\deg(G_\alpha).
$$
Remark that, if $s_*$ and $s'_*$ are the parabolic degree functions of 
quotient bundles $V$ and $W$ 
of $E_*$ of same rank $r''$, for all $\beta \in \R$ it is 
$s_{1+\beta}=s_\beta - r'' |I|$ and the same holds for $s'_*$, hence
the function $s_*-s'_*$ is periodic of period $1$ and 
$$
\int_\beta^{\beta +1}s_\alpha - s'_\alpha d\alpha = \deg(V_*) - \deg(W_*).\\
$$

\bigskip

\begin{definition} We define the scheme of parabolic quotients of $E_*$ of rank 
$r''$ and induced parabolic type $s_*$ as the stratum corresponding 
to the parabolic degree function $s_*$ and will denote it by $\pquot(E_*)$.\\
\end{definition}

\bigskip

\begin{remark}\label{par-qpar} \rm This construction allows to consider an 
algebraic structure on a scheme parametrizing parabolic quotient bundles, with 
possibly different underlying quasi-parabolic structures.\\

For instance, consider a parabolic bundle $E_*$ on $X$ at $I=\{p,q\}$
of rank $r$ and weight $0 < \alpha_1(p)=\alpha_1(q)=\alpha < 1$ and
suppose its underlying quasi-parabolic structure is such that 
$n_i(p), n_i(q)\ge r''$, for $i=1,2$. Let $s_*$ be a parabolic degree function
of quotient bundles of $E$ of rank $r''$ and degree $d''$, defined by
$s_0=d''$ and $s_\alpha=d''-r''$. \\

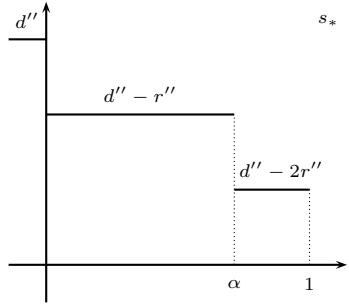
\begin{center}
\scriptsize
\psset{unit=0.5}
\pspicture(-1,-1)(2,7)
\psline{->}(-1,0)(8,0)
\psline{->}(0,-1)(0,7)
\psline(-1,6)(0,6)
\psline(0,4)(5,4)
\psline(5,2)(7,2)
\psline[linestyle=dotted,dotsep=1pt,linewidth=.5pt](5,0)(5,4)
\psline[linestyle=dotted,dotsep=1pt,linewidth=.5pt](7,0)(7,2)
\rput(-.5,6.5){$d''$}
\rput(2.5,4.5){$d''-r''$}
\rput(6.25,2.5){$d''-2r''$}
\rput(5,-.5){$\alpha$}
\rput(7,-.5){$1$}
\rput(7.5,6.5){$s_*$}
\endpspicture
\end{center}

Then the function $s_*$ corresponds to $r''+1$ quasi-parabolic structures, 
\ie as a set $\pquot(E_*)$ is the following disjoint union 
$$
\coprod_{k=0}^{r''}{\rm Quot}^{qpar}_{(k,r''-k),d''}(E,f_p)\cap
{\rm Quot}^{qpar}_{(r''-k,k),d''}(E,f_q).\\
$$
\end{remark}
\bigskip

Let $d_{r'',*}$ denote the minimal parabolic degree of a rank $r''$ 
quotient bundle of $E$. We denote by $\bar s_*$ a parabolic degree 
function whose parabolic degree is $d_{r'',*}$, that is 
$\int_{-1}^0\bar s_\alpha
d\alpha =d_{r'',*}$.\\

For a parabolic structure $*$ on $E$, we denote by $*'$ a parabolic structure 
obtained from $*$ by dropping one weight. 
For instance, consider the parabolic
structure obtained by dropping the highest weight $\alpha_L$:
$$
E \supset E_{\alpha_1}\supset E_{\alpha_2}\supset \dots \supset  
E_{\alpha_{L-1}} \supset E(-\sum_{p\in I}p).
$$
Then a parabolic degree function $s_*$ for quotient bundles of rank $r''$
completely determines the parabolic degree function for the structure $*'$
and we denote it by
$s'_{*'}$. This means that the stratum $\pquot(E_*)$  
naturally is a substratum of
the scheme ${\rm Quot}^{par}_{r'',s'_{*'}}(E_{*'})$.
Remark that the parabolic degree of these strata are such that 
\begin{eqnarray}\label{*to*}
\int_{-1}^0 s'_\alpha d\alpha=\int_{-1}^0 s_\alpha d\alpha -(\alpha_L-\alpha_{L-1})n''_{\alpha_{L+1}}.
\end{eqnarray}\\

\bigskip

\begin{remark}\label{max} \rm Consider a parabolic degree function 
$s_*$ for which
the last multiplicity is maximal, that is
$n''_{\alpha_{L+1}}=\min\{n_{\alpha_{L+1}},r''|I|\}$. Let $d''_*$
denote the parabolic degree of this stratum and $d''_{*'}$ denote
the parabolic degree of the stratum $s'_{*'}$. Then equality
(\ref{*to*}) is
$$
d''_{*'}=d''_* -(\alpha_L-\alpha_{L-1})n''_{\alpha_{L+1}}.
$$
Let $E \epi{} V$ be a quotient vector bundle of $E$ of rank $r''$.
Denote by $\bar d_*=\deg(V_*)$ and $\bar d_{*'}=\deg(V_{*'})$
the parabolic degrees with respect to the structures induced by
$*$ and $*'$ respectively. Then there exists an integer $n$ 
for which equality (\ref{*to*}) can be written as
$$
\bar d_{*'}=\bar d_{*}-(\alpha_L-\alpha_{L-1})n.
$$
Note that, since we assume that $n''_{\alpha_{L+1}}$ is maximal,
it is $n''_{\alpha_{L+1}} -n \le 0$.
This translates into the following inequality, on the differences of
parabolic degrees
\begin{eqnarray}\label{pardegdiff}
d''_{*'}-\bar d_{*'} = d''_{*}-\bar d_{*} - 
(\alpha_L-\alpha_{L-1})(n''_{\alpha_{L+1}}-n) \le d''_*-\bar d_*.
\end{eqnarray}
In particular, if $V$ is in a stratum corresponding to the minimal 
parabolic degree $d_{r'',*'}$, inequality (\ref{pardegdiff}) and
the fact that for all quotient bundle $V$ it is 
$\deg(V_*) \le d_{r'',*}$, imply the following inequality
$$
d''_{*'}-d_{r'',*'}=d''_{*'}-\deg(V_{*'}) 
\le d''_* - \deg(V_*) \le d''_* - d_{r'',*}.\\
$$
\end{remark}

\medskip

We still denote by $d''_*$ the parabolic degree of the stratum $s_*$.
We want to prove the following estimate for the dimension of the 
parabolic strata.\\

\bigskip

\begin{theorem} \label{cong} 
With the notations above, it is
$$
\dim(\pquot(E_*)) \le r''(r-r'') +r(d''_*-d_{r'',*}).\\
$$
\end{theorem}

\bigskip

\begin{remark} \rm This estimate depends on the parabolic invariants of the
stratum. In Appendix B we study the example of rank $2$ parabolic
bundles and show how, under some hypotheses, it is possible to 
get an estimate depending on the invariants
of the underlying vector bundles of the stratum.\\
\end{remark}

\bigskip

\pf The proof goes by induction on the number of weights. For $L=0$ the 
statement is given by the estimate of Popa and Roth 
\cite{popa-roth}, theorem 4.1 on the dimension
of Grothendieck's scheme of quotients. \\

We have to prove that the statement
holds for $L$ weights, provided it holds for $L-1$ weights. We prove it in 
two steps. The first one consists in proving the statement for $L$ weights
and under the assumption that $n''_{\alpha_{L+1}}$ is maximal. The second step
consists in drawing the general case from the first step.\\

\medskip
{\it First step}\\
Consider the parabolic structure $*'$ obtained from $*$ by dropping $\alpha_L$.
We still denote by $s'_{*'}$ the parabolic degree function
induced by $s_*$ with respect to the structure $*'$.
Then there is an obvious inequality
$$
\dim(\pquot(E_*)) \le \dim({\rm Quot}^{par}_{r'',s'_{*'}}(E_{*'}))
$$
and since the parabolic structure $*'$ has $L-1$ weights, the statement 
for the estimate of the dimension of ${\rm Quot}^{par}_{r'',s'_{*'}}(E_{*'})$
holds by the induction hypothesis. 
From this inequality and remark \ref{max} it follows
\begin{eqnarray*}
\lefteqn{\dim(\pquot(E_*))}\\ &\le 
\dim({\rm Quot}^{par}_{r'',s'_{*'}}(E_{*'}))
\le r''(r-r'') +r(d''_{*'}-d_{r'',*'})\\
& \le r''(r-r'') +r(d''_* - d_{r'',*}).
\end{eqnarray*}
This proves the statement for the strata of a parabolic structure with 
$L$ weights, whose last multiplicity is maximal.\\

\begin{remark} \label{trans} \rm We have actually proved that 
the statement holds for 
all strata for which there is an $i \in \{1,\dots,L+1\}$ 
such that $n''_{\alpha_i}$ is maximal. This is due to the fact that the
parabolic strata are obtained by successive flattening stratifications
and do not depend on the ``origin'' chosen for the filtration of the
parabolic bundle.
This translates into an isomorphism
$$
\pquot(E_*) \cong {\rm Quot}^{par}_{r'',s[\delta]_*}(E[\delta]_*)
$$
for all $\delta \in \R$,
where $E[\delta]_*$ is the parabolic bundle $E_*$ shifted by $\delta$
(see \cite{yokogawa}, definition 1.1) 
defined by $(E[\delta]_*)_\alpha=E_{\alpha +\delta}$ and $s[\delta]_*$
is the shifted parabolic degree function, that is 
$s[\delta]_\alpha=s_{\delta+\alpha}$.\\
\end{remark}

\medskip

{\it Second step}\\
We have to show that the result still holds when the multiplicities 
$n''_{\alpha_i}$ are such that no one of them is maximal.  
Of course, since $r''>0$ there is a nonzero multiplicity.
Without loss of generality, we can assume that it is
$$
0<n''_{\alpha_1}<\min\{n_{\alpha_1},r''|I|\}\le n_{\alpha_1}.
$$
We are going to need to add a (harmless) vector bundle
in the filtration of $E_*$, in order to 
conclude for the second step of the proof. \\

\bigskip

\begin{lemma}\label{locgen} There exists a vector bundle $\tilde E$ on $X$ 
such that
$E \supset \tilde E \supset E_{\alpha_1}$ and 
$\deg(\tilde E)=\deg(E_{\alpha_1})+n''_{\alpha_1}$ for which, if we consider
the parabolic structure $\tilde *$ obtained
from $*$ by adding $\tilde E$ to the structure $*$ with weight 
$\tilde \alpha \in ]0,\alpha_1[$, that is if $\tilde *$ is given by
$$
E \supset {\tilde E}=E_{\tilde \alpha}
\supset E_{\alpha_1}\supset E_{\alpha_2}\supset \dots \supset  E_{\alpha_L}
 \supset E(-\sum_{p\in I}p),
$$
then it is
$$
\dim(\pquot(E_*))=\dim({\rm Quot}^{par}_{r'',\tilde s_{\tilde *}}
(E_{\tilde *})),
$$
where $\tilde s_{\tilde *}$ is the parabolic degree function with same values
as $s_*$ at all $\alpha \in \{\alpha_0,\dots ,\alpha_L\}$
and $\tilde s_{\tilde \alpha}=d''=\tilde s _0$.\\
\end{lemma}

\bigskip

This lemma allows to finish the induction argument. Fix a parabolic 
structure $\tilde *$ as in the lemma and consider $\tilde * '$ the parabolic 
structure obtained from $\tilde *$ by dropping $\alpha_0$. The parabolic 
structure $\tilde * '$ has $L$ weights and the stratum associated with
$\tilde s '_{\tilde *'}$ has
maximal first multiplicity. Remark 
that the parabolic degree $d''_{\tilde * '}$
of the stratum $\tilde s '_{\tilde *'}$
is such that
$$
d''_{\tilde * '}=d''_{\tilde *}=d''_* + \tilde \alpha n''_{\alpha_1}.\\
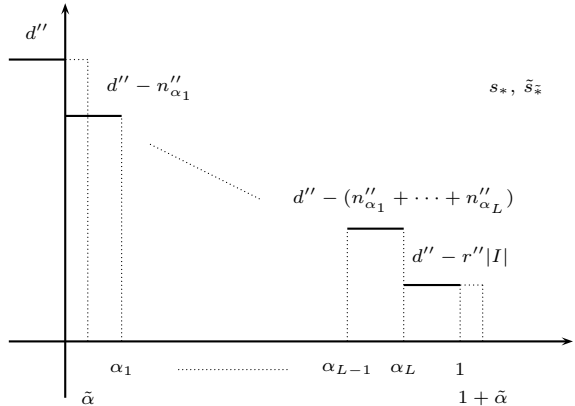
$$

\begin{center}
\scriptsize
\psset{unit=0.75}
\pspicture(-1,-1)(9,6)
\psline{->}(-1,0)(9,0)
\psline{->}(0,-1)(0,6)
\psline(-1,5)(0,5)
\psline[linestyle=dotted,dotsep=1pt,linewidth=.5pt](0,5)(.4,5)
\psline(0,4)(1,4)
\psline[linestyle=dotted,dotsep=1pt,linewidth=.5pt](1.5,3.5)(3.5,2.5)
\psline(5,2)(6,2)
\psline(6,1)(7,1)
\psline[linestyle=dotted,dotsep=1pt,linewidth=.5pt](.4,0)(.4,5)
\psline[linestyle=dotted,dotsep=1pt,linewidth=.5pt](1,0)(1,4)
\psline[linestyle=dotted,dotsep=1pt,linewidth=.5pt](5,0)(5,2)
\psline[linestyle=dotted,dotsep=1pt,linewidth=.5pt](6,0)(6,2)
\psline[linestyle=dotted,dotsep=1pt,linewidth=.5pt](7,0)(7,1)
\psline[linestyle=dotted,dotsep=1pt,linewidth=.5pt](7,1)(7.4,1)
\psline[linestyle=dotted,dotsep=1pt,linewidth=.5pt](7.4,0)(7.4,1)
\rput(-.5,5.5){$d''$}
\rput(1.5,4.5){$d''-n''_{\alpha_1}$}
\rput(6,2.5){$d''-(n''_{\alpha_1}+\dots +n''_{\alpha_L})$}
\rput(7,1.5){$d''-r''|I|$}
\rput(.4,-1){$\tilde\alpha$}
\rput(1,-.5){$\alpha_1$}
\psline[linestyle=dotted,dotsep=1pt,linewidth=.5pt](2,-.5)(4,-.5)
\rput(5,-.5){$\alpha_{L-1}$}
\rput(6,-.5){$\alpha_L$}
\rput(7,-.5){$1$}
\rput(7.4,-1){$1+\tilde\alpha$}
\rput(8,4.5){$s_*, \, \tilde s_{\tilde *}$}
\endpspicture
\end{center}

We have the inclusion 
$$
{\rm Quot}^{par}_{r'',\tilde s_{\tilde *}}(E_{\tilde *})
\subseteq 
{\rm Quot}^{par}_{r'',\tilde s'_{\tilde *'}}(\tilde E_{\tilde *'}).
$$ 
From this inclusion and lemma \ref{locgen} we get the inequality
$$
\dim(\pquot (E_*))=\dim({\rm Quot}^{par}_{r'',\tilde s_{\tilde *}}
(E_{\tilde *}))
\le \dim({\rm Quot}^{par}_{r'',\tilde s'_{\tilde *'}}
(\tilde E_{\tilde *'})).\\
$$

By the first step of the induction argument, the last dimension is less than 
or equal to 
$$
r''(r-r'')+r(d''_{\tilde *'}-d_{r'',\tilde *'}). 
$$
Let 
$\bar s_{\tilde * '}$ be a parabolic degree function realizing the
minimal parabolic degree $d_{r'',\tilde *'}$. By remark \ref{max}, it is
$$
d''_{\tilde *'}-d_{r'',\tilde *'} \le d''_{\tilde *}-\bar d_{\tilde *},
$$ 
where $\bar d_{\tilde *}$ is the parabolic degree of some vector bundle 
$\bar V$ of the stratum $\bar s_{\tilde * '}$. Moreover, if we denote by 
$\bar d_*$ the parabolic degree of $\bar V$ with respect to 
the original parabolic structure $*$, since the structure
$\tilde *$ has one weight more than $*$, it is 
$\bar d_{\tilde *}\ge \bar d_*$. Thus we get
$$
d''_{\tilde *'}-d_{r'',\tilde *'} \le d''_{\tilde *}-\bar d_{\tilde *}
\le
d''_{\tilde *}-\bar d_*.
$$
Recalling the expression of $d''_{\tilde*}$, we get the inequality
$$
\dim(\pquot(E_*)) \le r''(r-r'') + 
r(d''_*+\tilde \alpha n''_{\alpha_1} -\bar d_*)
$$
$$
\le r''(r-r'') +r(d''_* - d_{r'',*} +\tilde \alpha n''_{\alpha_1}).
$$
Then in order to get the inequality of the statement 
it is enough to remark that, according to lemma \ref{locgen}, we can choose
$\tilde \alpha$ as small as we like. \cvd\\

\bigskip

{\noindent \emph{\mbox{Proof of lemma \ref{locgen}:\ \ \ }}}
It will be enough to prove that
$$
\dim(\pquot(E_*))\le
\dim({\rm Quot}^{par}_{r'',\tilde s_{\tilde *}}(E_{\tilde *})),
$$
that is, to find some subscheme of the substratum associated with 
$\tilde s_{\tilde *}$ with same dimension as the stratum associated with
$s_*$. Let $Q$ be an irreducible component of $\pquot(E_*)$ of
maximal dimension and let $V$ be a quotient bundle of $E$, 
representing a point of $Q$. By assumption on the function $s_*$ it is
\begin{eqnarray*}
\deg(E_{\alpha_1}/E(-\sum_{p\in I}p)) = n_{\alpha_2} + \dots 
+ n_{\alpha_{L+1}} > \\
\ \ \deg(V_{\alpha_1}/V(-\sum_{p\in I}p))=n''_{\alpha_2} 
+ \dots + n''_{\alpha_{L+1}}
\end{eqnarray*}
and moreover, since $n''_{\alpha_1}>0$, it is $n''_{\alpha _2} + \dots + 
n''_{\alpha_{L+1}} < r''|I|$.
Translating this into quasi-parabolic conditions, this means that there 
are some $p\in I$ such that $\alpha_1(p)=\alpha_1$ and at $p$ the
induced parabolic structure is a strict inclusion in the fiber of $V$:
$$\xymatrix{
E_p \ar@{->>}^{q}[r] & V_p \\
E_{p,2} \ar@{^{(}->}^{i}[u] \ar@{->>}[r] & V_{p,2}=\im(qi)\ar@{^{(}->}[u]. \\
}$$
On the other hand, we can add the missing generators of $V_p$ at such points:
we can choose a subset $I' \subseteq I$ of points such that 
$\alpha_1(p)=\alpha_1$ and at these points add a linear subspace $H_{p,2}$
to $E_{p,2}$, with $\rk(H_{p,2})=\rk(V_p/V_{p,2})$ for which the composition
$$
H_{p,2}\oplus E_{p,2} \mono{} E_p \epi{} V_p
$$
has rank $r''$ and $\sum_{p\in I'}\rk(H_{p,2})=n''_{\alpha_1}$. 
It is enough to choose the vector space $H_{p,2}$ as the image of a section
of the surjective linear map
$$
E_p \epi{} V_p/V_{p,2}.
$$
The parabolic structure $\tilde *$ is obtained by enriching
the flags of the quasi-parabolic structures at the points of $I'$ 
in the following way:
$$
E_p \supset H_{p,2}\oplus E_{p,2} \supset E_{p,2} \supset
\dots \supset E_{p,l_p} \supset 0
$$
with weights 
$(\alpha_i(p))=(\tilde \alpha , \alpha_1(p),\dots,\alpha_{l_p}(p))$. 
This means that $V$ represents a point of the substratum
${\rm Quot}^{par}_{r'',\tilde s_{\tilde *}}(E_{\tilde *})$. All is left to
check is that this is true in an open neighbourhood of $V$ in $Q$.\\

By semi-continuity of the rank, the composed linear map 
$$
H_{p,2}\oplus E_{p,2} \mono{} E_p \epi{} W
$$ 
has rank $r''$, for all vector spaces $W$ in an open neighbourhood 
$U_p$ of the isomorphism class of the
fiber $V_p$ in the grassmannian $G_p=\grass_{r''}(E_p)$. Hence the 
induced quasi-parabolic 
filtration is of the same type, for all vector bundle in an open 
neighbourhood of $V$. Recall that all the points of the parabolic strata 
are vector bundles and consider, for all
$p \in I'$, the map
$$
\epsilon_p: Q \ra G_p, \ \ \ [E \epi{\pi} W] \mapsto [E_p \epi{\pi_p} W_p].
$$
Let $Q'$ be the open subscheme of $Q$ defined as the intersection
$$
Q'=\bigcap_{p\in I'}\epsilon_p^{-1}(U_p).
$$
Since $Q$ is irreducible, we have $\dim(Q')=\dim(Q)$ and by construction 
$Q'$ is a subscheme of 
${\rm Quot}^{par}_{r'',\tilde s_{\tilde *}}(E_{\tilde *})$. \cvd \\

\bigskip

\section{\bf {Sections of the line bundle $\cl ^{par}$}}

\bigskip

Let $\cm^{par}$ denote the moduli space of parabolic bundles on $X$ of
rank $r$, trivial determinant and parabolic structure at $I$ of
multiplicities $((n_1(p), \dots , n_{l_p}(p))_{p\in I})$ and weights
$$
0 \le \alpha_1(p)=0 < \alpha_2(p) < \alpha_3(p)< \dots <
\alpha_{l_p}(p) < 1.
$$
Let $((d_1(p), \dots, d_{l_p-1}(p))_{p\in I},k)$ be strictly positive
integers such that, for all $j=2, \dots, l_p$, the weights at $p$ can
be written as
$\alpha_j(p)=\displaystyle\frac{1}{k}\sum_{h=1}^{j-1}d_h(p)$.\\

A family $\ce_*$ of parabolic bundles at $I$, of rank $r$, trivial
determinant, multiplicities $((n_i(p))_{p\in I})$ and weights
$((d_j(p))_{p\in I}, k)$ parametrized by a scheme $S$ is a vector
bundle $\ce$ over $X\times S$ of rank $r$, such that $\det(\ce)=\co_{X\times S}$ and,
for all $p\in I$, quotient bundles $Q_i(p)$ of $\ce _{|\{p\}\times
S}$, of rank $r_i(p)=n_1(p)+\dots +n_i(p)$, such that, by letting $\ck
_i(p)$ denote the kernel
$$\xymatrix{ 0 \ar[r] & \ck _i (p)= \ker (\pi_i(p)) \ar[r] & \ce
_{|\{p\}\times S} \ar[r]^{\pi_i(p)}  & Q_i(p) \ar[r] & 0,\\ }$$ then
$\ck _i(p) \subset \ck _{i-1}(p)$ for all $i=1 \dots , l_p-1$. The
family is parabolic in the sense that, for all $s \in S$ the vector
bundle $\ce _s$, has the quasi-parabolic structure
$$
{\ce _s}_p \supset {\ck _1}(p)_s \supset \dots \supset {\ck
_{l-1}}(p)_s \supset 0
$$ 
and weights $((d_j(p))_{p\in I},k)$. Actually the family $\ce$ has a
weighted filtration, induced by the quotients $Q_i(p)$, that is
$$
\ce = \ce _0 \supset \ce _{\delta_1(p)} \supset \dots \supset \ce
_{\delta_1(p)+ \dots + \delta_{l_p-1}(p)} \supset \ce _1 = \ce \ts
\pi_X^*\co_X(-p),
$$
where $\delta_j(p)=\displaystyle\frac{d_j(p)}{k}$ and $Q_j(p) \cong \ce
/\ce _{\delta_1(p)+ \dots + \delta_j(p)}$. Suppose that the family is
semistable and let $\phi_S : S \ra \cm^{par}$ the modular
morphism. Suppose that $\displaystyle\sum_{p\in
I}\displaystyle\sum_{i>j} \displaystyle\frac{n_i(p)d_j(p)}{r}\in \Z$
and let $\cl ^{par}(\ce_*)$ be the line bundle on $S$ defined as the tensor product
$$
\cl ^{par}(\ce_*)=(\det R\pi_S\ce)^{\ts k}\ts \Ts_{p\in I}\Ts_j (\det
Q_j(p))^{\ts d_j(p)}\ts (\det \ce _{|\{q\}\times S})^{\ts e}.
$$
Here $e$ is an integer depending on the parabolic structure
defined by
$$
e=\frac{1}{r}\sum_{p\in I}\sum_{i>j}n_i(p)d_j(p)+k(1-g)-\sum_{p\in
I}\sum_jd_j(p).\\
$$

\bigskip

\begin{theorem} (\cite{narasimhan-ramadas}, theorem 1, for rank $2$; \cite{pauly 1},
  theorem 3.3, for arbitrary rank) There exists a unique ample line bundle $\cl ^{par}$
  over $\cm ^{par}$ such that, for all semistable parabolic family
  $\ce_*$ parametrized by $S$, it is $\phi_S^*\cl ^{par}=\cl
  ^{par}(\ce_*)$.\\
\end{theorem}

\bigskip

In the case of rank two parabolic bundles, Pauly gives in \cite{pauly
  2} a method to produce sections of $\cl ^{par}$ of type theta. In
  what follows, we extend his method to the rank $r$ case and produce
  sections of  ${\cl ^{par}}^{\ts h}$, for all $h \in \N$. \\

Let $\ce _*$, $\cf _*$ be families of parabolic bundles on $I$,
parametrized by $S$, of quotients and flags respectively
\[
\begin{array}{cc}
0 \arr {} \ce _\alpha \arr {i_\alpha} \ce \arr {p_\alpha} Q_\alpha
\arr {} 0 & 0 \ra \ck_i(p) \ra \ce _{|\{p\}\times S} \ra Q_i(p) \ra 0,
\\ 0 \arr {} \cf _\alpha \arr {i'_\alpha} \cf \arr {p'_\alpha}
Q'_\alpha \arr {} 0 &  0 \ra \ck'_j(p) \ra \cf _{|\{p\}\times S} \ra
Q'_j(p) \ra 0.
\end{array}
\]
A morphism $\varphi : \ce \ra \cf$ of vector bundles is parabolic if
the composed map $p'_\alpha \varphi i_\alpha$
is the zero morphism, for all $\alpha \in \R^+$. 
The sheaf of parabolic homomorphisms
is a locally free subsheaf of $\fhom(\ce,\cf)$, that will be denoted by
$\fhom(\ce_*,\cf_*)$. The quotient sheaf
$\coker(\fhom(\ce_*,\cf_*) \mono{} \fhom(\ce,\cf))$ 
is a family of torsion sheaves
parametrized by $S$ whose support is contained in the parabolic subset $I$,
that we denote by $T_{\ce_*,\cf_*}$.\\ 

The sheaf $\fhom (\ce _*, \cf_*)$ is the family on $S$
parametrizing the sheaves of parabolic morphisms between bundles of
the families, that is for all $s \in S$ there is a natural isomorphism
$$
\fhom({\ce_*}_s, {\cf_*}_s) \cong \fhom(\ce_*,\cf_*)_s.\\
$$

Let $F$ be a vector bundle on $X$ of rank $hk$ such that 
$$
\deg(F)=\frac{h}{r}\sum_{p\in I}\sum_{i>j}n_i(p)d_j(p) +
hk(g-1)
$$ 
and let $F_*$ be a parabolic structure at $I$ of
multiplicities 
$$
((hd_1(p), \dots, hd_{l_p-1}(p), hk - h \sum_jd_j(p))_{p\in
I})
$$ 
and weights $((d_j(p))_{p\in I}, k)$. Let $\pi_X^*F_*$ be the
constant family of parabolic bundles of value $F_*$, parametrized by
$S$. For $s\in S$, let ${\ce_*}_s$ denote the parabolic bundle of $\ce_*$
over $s$. Then it is
\begin{eqnarray*}
\lefteqn{\chi({\ce_*}_s,F_*)}\\  & = rhk(\mu({\ce_*}_s)+(g-1))
-hk\deg(\ce_s)+rhk(1-g) -  \displaystyle\sum_{p\in
I}\displaystyle\sum_{i>j}n_i^{{\ce_*}_s}(p)n_j^{F_*}(p)\\ & = h
\displaystyle\sum_{p\in
I}\sum_{i>j}n_i(p)d_j(p)-\displaystyle\sum_{p\in
I}\displaystyle\sum_{i>j}n_i(p)hd_j(p)=0.
\end{eqnarray*}\\

Fix a basis of $F_p= \displaystyle\bigoplus_j F_{p,j}/F_{p,j+1}$ and
let $\ct_{\ce_*,F_*}$ denote the family of torsion sheaves of
the short exact sequence of parabolic morphisms
\begin{eqnarray}\label{morfpar}
0 \ra  \fhom(\ce _*,\pi_X^*F_*) \ra \fhom(\ce,\pi_X^*F) \arup{\rm p}
 \ct_{\ce_*,F_*} \ra 0.
\end{eqnarray}\\

\begin{lemma} With the notations above, it is
$$
\det R \pi_S\fhom(\ce_*,\pi_X^*F_*)\cong \cl ^{par}(\ce_*)^{\ts h}.
$$\\
\end{lemma}

\pf By the short exact sequence (\ref{morfpar}) there is a natural
isomorphism
$$
\det R \pi_S\fhom(\ce_*,\pi_X^*F_*)\cong (\det R
\pi_S\fhom(\ce,\pi_X^*F))\ts (\det R\pi_S\ct_{\ce_*,F_*})^\vee.
$$
By Serre duality theorem, \cite{pauly 2}, lemma 3.4, there is an
isomorphism
$$
\det R \pi_S\fhom(\ce,\pi_X^*F) = \det R \pi_S \ce ^\vee \ts \pi_X^*F
\cong \det R \pi_S \ce \ts \pi_X^*(F^\vee \ts {\rm K}_X).
$$
The vector bundle $\det\ce _{|\{q\}\times S}$ is independent of $q \in
X$ and by \cite{pauly 2}, lemma 3.5 it follows that
$$
\det R \pi_S \ce \ts \pi_X^*(F^\vee \ts {\rm K}_X) \cong (\det R \pi_S
\ce)^{\ts hk} \ts (\det\ce _{|\{q\}\times S})^{\ts -\deg (F^\vee\ts
{\rm K}_X)}.
$$ 
Now, since the degree of $F^\vee\ts {\rm K}_X$ can be computed as 
$$
\deg (F^\vee\ts {\rm K}_X)= hk \deg ({\rm K}_X) -
\deg(F)= hk(g-1)-\displaystyle\frac{h}{r}\sum_{p \in
I}\sum_{i>j}n_i(p)d_j(p),
$$ 
the first determinant bundle is isomorphic to
\begin{eqnarray*}
\lefteqn{\det R \pi_S\fhom(\ce,\pi_X^*F)\cong(\det R \pi_S \ce)^{\ts hk}} \\ & 
\ts (\det\ce _{|\{q\}\times S})^{\ts \frac{h}{r}\sum_{p\in
I}\sum_{i>j}n_i(p)d_j(p) +hk(1-g)}.
\end{eqnarray*}\\

The sheaf $\ct_{\ce_*,F_*}$ is a family of skyscraper sheaves
supported at $I$, hence the sheaf $R^1{\pi_S}_*\ct_{\ce_*,F_*}$ is
zero and there is an isomorphism
$$
{\pi_S}_*\ct_{\ce_*,F_*}\cong \bigoplus_{p\in I}\bigoplus
_{j=1}^{l_p-1}{\ck_j(p)}^\vee \ts \co_S ^{hd_j(p)}.
$$
Thus the second determinant can be computed as follows
\begin{align*}
(\det R \pi_S \ct_{\ce_*,F_*})^\vee & \cong \det\bigoplus_{p\in
I}\bigoplus _{j=1}^{l_p-1}{\ck_j(p)}^\vee \ts \co_S^{hd_j(p)} \cong
\Ts_{p\in I}\Ts _{j=1}^{l_p-1}\det ({\ck_j(p)}^\vee\ts
\co_S^{hd_j(p)})\\ & \cong \Ts_{p\in I}\Ts
_{j=1}^{l_p-1}(\det{\ck_j(p)}^\vee)^{\ts {hd_j(p)}}.
\end{align*}\\
By definition, it is $\ck_i(p)=\ker(\pi_i(p))$ and this yields the
isomorphism
$$
\Ts_{p\in I}\Ts _{j=1}^{l_p-1}(\det{\ck_j(p)}^\vee)^{\ts
{hd_j(p)}}\cong \Ts_{p\in I}\Ts _{j=1}^{l_p-1}(\det Q_j(p))^{\ts
{hd_j(p)}}\ts (\det\ce ^\vee_{|\{p\}\times S})^{\ts {hd_j(p)}}.
$$\\
The lemma then follows from the fact that for all
$p\in I$ there is a natural isomorphism
$\det\ce_{|\{q\}\times S}\cong \det\ce_{|\{p\}\times S}$ and the equality
$$
he=\frac{h}{r}\sum_{p\in I}\sum_{i>j}n_i(p)d_j(p) +hk(1-g)-h\sum_{p\in
I}\sum_jd_j(p).\\
$$
\cvd\\

Let $\cf$ be a quasi-coherent $\co _{X \times S}$-module, flat over
$S$.  Recall how one obtains a complex that is quasi-isomorphic to $R
\pi_S \cf$. By the relative version of Serre A theorem, there is an
integer $m_0$ such that, if $m\ge m_0$, the natural evaluation morphism
$$
K_0=\pi_X^*\co_X(-m)\ts \pi_X^* {\pi_X}_*\cf(m)\epi q \cf
$$
is surjective. Since $\deg(\co_X(-m))<0$, it is ${\pi_S}_*K_0=0$ and
if we denote by $K_1=\ker q$ it is ${\pi_S}_*K_1=0$ as well. Moreover the
higher direct image sheaves $\cl_1=R^1{\pi_S}_*K_1$,
$\cl_0=R^1{\pi_S}_*K_0$ are locally free. Hence the short exact
sequence
$$
\xymatrix{ 0 \ar[r] & K_1 \ar[r]^{i} & K_0 \ar[r]^{q} & \cf \ar[r] &
0 \\   
}$$ 
yields the long exact sequence in cohomology
$$\xymatrix{ 
0 \ar[r] & {\pi_S}_*\cf \ar[r] & R^1{\pi_S}_*K_1
\ar[r]^{R^1{\pi_S}_*i} & R^1{\pi_S}_*K_0 \ar[r] & R^1{\pi_S}_*\cf
\ar[r] & 0 \\ 
}$$ 
and there is a natural isomorphism $\det R\pi_S \cf
\cong \det \cl _0 \ts \det \cl _1^\vee$.\\

Let $0 \ra \cl _1 \arup \nu \cl _0 \ra 0$ be a complex of locally free
sheaves on $X\times S$, quasi-isomorphic to $R
\pi_S\fhom(\ce_*,\pi_X^*F_*)$. The hypothesis on the Euler
characteristic $\chi({\ce_*}_s,F_*)=0$, for all $s \in S$, is equivalent 
to the assumption that the locally free sheaves
$\cl_i$ have same rank and the morphism of
vector bundles $\det \nu : \det \cl _1 \ra \det \cl _0$ defines a
section of $(\det\cl _1)^\vee \ts \det \cl _0 = \det R
\pi_S\fhom(\ce_*,\pi_X^*F_*)\cong \cl ^{par}(\ce_*)^{\ts h}$, that we
denote by $\theta_{F_*}^{\ce _*}$. This section is zero at a point $s
\in S$, if and only if
$$
\dim (\Hom({\ce_*}_s,F_*))=\dim (\ext({\ce_*}_s,F_*)) \neq 0. \\
$$

To show that this produces a section of the line bundle on the moduli
space $\cm ^{par}$, recall its construction (see, for instance,
\cite{pauly 1}, theorem 2.3). Let $Q$ be the scheme of quotients of
rank $r$ and trivial determinant of $\co_X^{P(n)}(-n)$, where $P$ is
the Hilbert polynomial of such quotients and $n$ is an integer, $n \gg
0$. Let $\Omega$ denote the open subset of $Q$ of locally free
quotients, $\cf$ the universal family of quotients on $X\times \Omega$
and ${\rm F}_p$ the flag varieties bundle of multiplicities $(n_1(p),
\dots , n_{l_p}(p))$
$$
{\rm F}_p = \Flag _{(n_1(p), \dots , n_{l_p}(p))}(\cf _{|\{p\} \times
\Omega}) \epi {\pi(p)} \Omega.
$$
Let $\cq_i(p)$ denote the universal quotients on ${\rm F}_p$ and let
${\mathcal R}$ be the fibred product of the ${\rm F}_p$'s, for $p \in
I$, over $\Omega$. We still denote by $\cf_*$ and $\cq_i(p)$ the
universal families obtained by pullback to ${\mathcal R}$. The
parabolic family $\cf_*$, with parabolic quotients $(\cq_i(p))$, is
locally a universal family of parabolic bundles. Let ${\mathcal
R}^{ss}$ be the open subscheme of ${\mathcal R}$ of semistable
parabolic bundles. Then $\cm ^{par}$ is obtained as the good quotient
of ${\mathcal R}^{ss}$ for the natural action of $SL(P(n))$.\\

Consider the line bundle $\cl ^{par}(\cf_*)$ on ${\mathcal
  R}^{ss}$. By \cite{pauly
  1}, theorem 3.3 it descends to the moduli space $\cm ^{par}$. 
The section $\theta_{F_*}^{\cf_*}$ is
  $SL(P(n))$-invariant, thus it descends to a section of ${\cl
  ^{par}}^{\ts h}$, that will be called {\it parabolic theta function}
  (of order $h$) associated with the parabolic bundle $F_*$. \\

\bigskip

\section{\bf{Zeroes of parabolic theta functions}}

\bigskip

Let $E_*$ be a semistable parabolic bundle on $X$ at $I$, of rank $r$,
trivial determinant, multiplicities $((n_i(p))_{p \in I})$ and weights
$((d_j(p))_{p \in I},k)$. For a parabolic bundle $F_*$ on $X$ at $I$,
of rank $\ell k$, slope $\mu(E_*)+g-1$, multiplicities $((\ell
d_1(p),\dots, \ell d_{l_p-1}(p),\ell
(k-\displaystyle\sum_{i=1}^{l_p-1}d_i(p)))_{p \in I})$ and same
weights as $E_*$, the parabolic theta function associated with $F_*$ is
zero at the point of $\cm^{par}$ represented by $E_*$, if and only if
$\Hom(E_*,F_*)={\rm H}^0(\fhom(E_*, F_*)) \neq \{0\}$. Let
$d_{l_p}(p) =k-\displaystyle\sum_{i=1}^{l_p-1}d_i(p)$ and let $\cm_\ell^{'par}$ denote the moduli space of equivalence classes of semistable parabolic bundles $F_*$, with which we can associate parabolic theta functions of order $\ell$. Recall that its dimension is given by
$$
\dim(\cm_\ell^{'par})=(\ell k)^2(g-1) + \sum_{p\in I}d_{\ell
  d_1(p),\dots,\ell d_{l_p}(p)}+1. 
$$
Let $r''$ be an integer such that $0<r''\le r$ and let $\ce _{r''}$
denote the family of isomorphism classes of stable parabolic bundles
$F_*$ such that there is a morphism $\varphi_* :E_* \ra F_*$ of rank
$r''$. We prove in this section that whenever 
$\ell \ge r''(r-r'')$ and $\ell \ge \frac{r}{k}$, 
then
$$
\dim(\ce _{r''})\le (\ell k)^2(g-1) + \sum_{p\in I}d_{\ell
  d_1(p),\dots,\ell d_{l_p}(p)}.
$$
This will prove theorem \ref{gg} since if $I \neq \emptyset$ 
then $k\ge 2$ and
$$
\sup_{0<r''\le r}\left\{r''(r-r''),\frac{r}{k}\right\}
\le \left[\frac{r^2}{4}\right].
$$
Then there exists a nonempty open subset $\cu$ of
the moduli space $M_\ell^{'par}$, such that for
all $F_*$ representing a point of $\cu$ it is $\Hom(E_*,F_*)= 0$.\\

\bigskip

\subsection{Images of parabolic morphisms} \hfill

\bigskip

Let $\varphi_* :E_* \ra F_*$ be a morphism. The image of $\varphi$ is
a quotient bundle of $E$, denote it by $V=\im(\varphi)$ and let $V_*$ be the
induced parabolic structure via the quotient morphism $E_* \epi {} V$.
The subbundle $V'$ of $F$ generated by $V$ inherits a natural
parabolic structure as well, via the injective morphism to $F_*$. We want to 
compare these two induced parabolic structures. Note that, if the support of the quotient sheaf $V'/V$ does not intersect the parabolic subset $I$, the two parabolic structures necessarily have the same multiplicities.\\
Suppose for simplicity that $I=\{p\}$ and
let $n'_i= \deg (V'_{\alpha_{i-1}}/V'_{\alpha_i})$ be the
multiplicities of the parabolic structure induced on $V'$ by $F_*$ and
$n''_i= \deg (V_{\alpha_{i-1}}/V_{\alpha_i})$ be the multiplicities
induced on $V$ by $E_*$.\\

\begin{proposition} With these notations, it is
\begin{eqnarray}\label{var-cell}
\lefteqn{\deg(V) + \frac{1}{k}\sum_j(r''-r''_j)d_j \le} \\ 
& \displaystyle  \deg(V) +
\frac{1}{k}\sum_j(r''-r'_j)d_j \le \deg(V') +
\frac{1}{k}\sum_j(r''-r'_j)d_j. \notag
\end{eqnarray}
In particular, the following inequality holds
$$
\deg(V_*) \le \deg(V'_*).
$$\\
\end{proposition}

\pf  We are actually going to prove that $\deg
(V/V_{\alpha_i})\ge\deg  (V'/V'_{\alpha_i})$, for all $i=1, \dots ,
l$. In fact, this inequality can be rewritten as
$$
r''_i=\displaystyle\sum_{j \le i}n''_j \ge
r'_i=\displaystyle\sum_{j\le i}n'_j, 
$$
for all $i$. The underlying
vector bundle $V'$ is the saturation of $V$ in $F$ and so $\deg(V) \le
\deg(V')$. From these facts inequalities (\ref{var-cell}) follow.\\

Let $i \in \{1, \dots, l\}$ and denote $\alpha_i=\alpha$.  The
morphism $\varphi$ is parabolic, so $V_\alpha \cong
\im(\varphi_\alpha)$ and the diagram (\ref{A})  commutes. From this we
deduce the commutative diagram (\ref{B}), hence a morphism $j_\alpha:
V_\alpha \ra V'\times_{F}F_\alpha\cong V'_\alpha$. \\

\begin{multicols}{2}
\begin{eqnarray}\label{A}\xymatrix{
E \ar[rr]^{\varphi} \ar@{->>}[dr] & & F\\ & V=\im(\varphi)
\ar@{^{(}->}[ur] &\\ E_\alpha \ar[rr]^{\varphi_\alpha}|\hole
\ar@{->>}[dr] \ar@{^{(}->}[uu] & & F_\alpha \ar@{^{(}->}[uu] \\  &
V_\alpha=\im(\varphi_\alpha) \ar@{^{(}->}[ur] \ar@{^{(}->}[uu]& \\ }
\end{eqnarray}
\begin{eqnarray}\label{B} \xymatrix{
F_\alpha \ar@{^{(}->}[rrr] & & & F \\ V_\alpha \ar@{^{(}->}[u]
\ar@{^{(}->}[r] & V\ar@{^{(}->}[dr] \ar@{^{(}->}[rr] & & F
\ar@{=}[u]\\ & & V'\ar@{^{(}->}[ur] &  \\ }
\end{eqnarray} 
\end{multicols}

The morphism $j_\alpha$ is such that $i_{V'}^\alpha j_\alpha$ is
injective, so it is injective as well. Thus the cokernel $\tau_\alpha$
has rank zero.  Denote by $G'_\alpha =\coker(i_\alpha)$ and $G_\alpha
=\coker(V'_\alpha \ra F_\alpha)$ and let $i_G^{\alpha}:G_\alpha \ra G$
be the morphism of the induced parabolic structure on $G$. From the
commutativity of the cube, we deduce the morphisms $v:\tau_\alpha \ra
\tau$ and $i_{G'}^{\alpha}:G'_\alpha \ra G'$ as well as $\tau_\alpha \ra
G'_\alpha$ this translates into the diagram (\ref{completo}).  \\

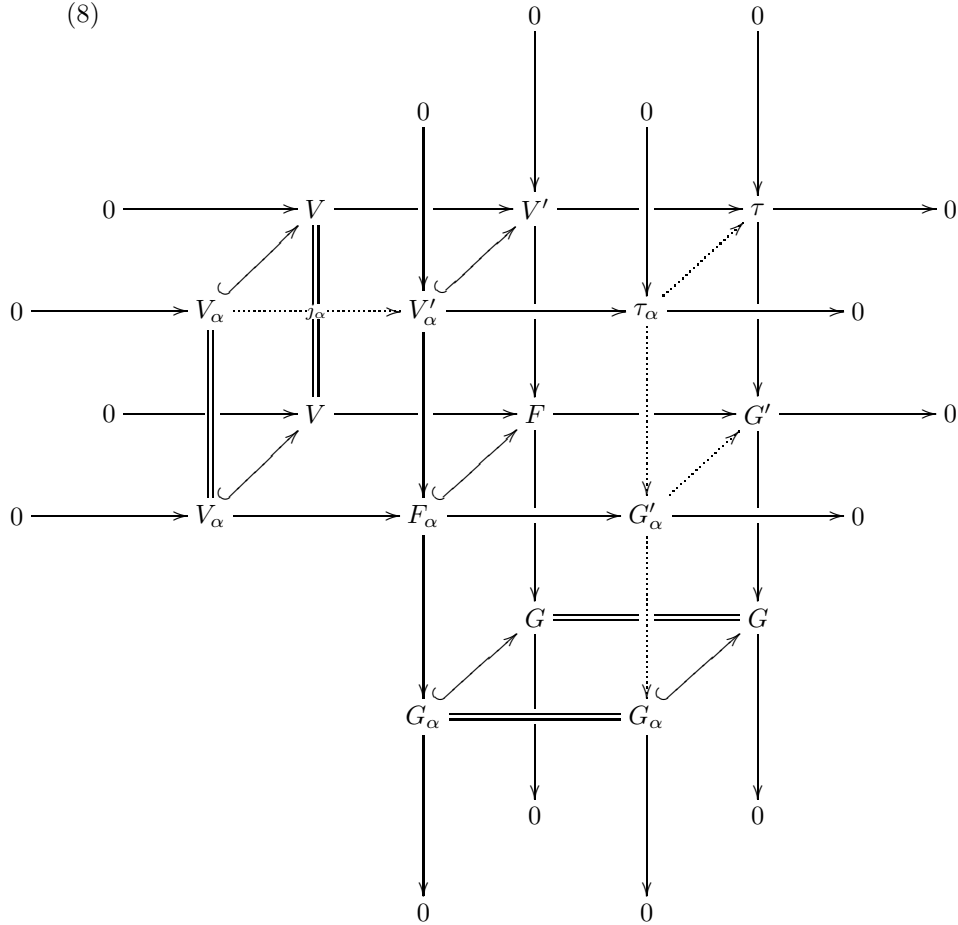
\begin{eqnarray}
\label{completo}\xymatrix{
& & & & & 0 \ar[dd] & & 0 \ar[dd] & & \\ & & & & 0 \ar[dd]&         &
0 \ar[dd] &     & & \\  & 0 \ar[rr] & & V \ar@{=}[dd]|\hole
\ar[rr]|\hole & & V' \ar[dd]|\hole \ar[rr]|\hole & & \tau
\ar[dd]|\hole \ar[rr] & & 0 \\ 0 \ar[rr] & & V_\alpha \ar@{=}[dd]
\ar@{.>}[rr]|{j_\alpha} \ar@{^{(}->}[ur] & &  V'_\alpha \ar[dd]
\ar[rr] \ar@{^{(}->}[ur] & & \tau_\alpha \ar@{.>}[dd]  \ar[rr]
\ar@{.>}[ur] & &0 &\\ & 0 \ar[rr]|\hole & & V \ar[rr]|\hole & & F
\ar[dd]|\hole \ar[rr]|\hole & & G' \ar[rr]\ar[dd]|\hole & & 0 \\ 0 \ar
[rr] & & V_\alpha \ar@{^{(}->}[ur]\ar[rr] & & F_\alpha \ar[rr] \ar[dd]
\ar@{^{(}->}[ur] & & G'_\alpha \ar@{.>}[ur]\ar[rr]\ar@{.>}[dd] & & 0 &
\\ & & & & & G \ar@{=}[rr]|\hole \ar[dd]|\hole & & G \ar[dd] &  \\ & &
& & G_\alpha \ar@{=}[rr] \ar@{^{(}->}[ur] \ar[dd] & & G_\alpha
\ar@{^{(}->}[ur]\ar[dd] & &  &\\ & & & & & 0 & & 0  & & \\  & & & & 0
& & 0  & & & \\  }
\end{eqnarray}\\

By the snake lemma it
follows that the morphism $\tau_\alpha \ra G'_\alpha$ is injective and
its cokernel is isomorphic to $G_\alpha$. These morphisms are such
that in the diagram (\ref{completo}) each horizontal and vertical
diagram is commutative.
Starting over this process
from the vertical diagram of weight $\alpha$ thus
obtained, we can add the corresponding vetical diagram of weight
$1$.  Now, the first nontrivial horizontal diagram is (\ref{ff}) 
and the morphism $\tau \arup{u} \tau_\alpha \arup{v} \tau$ is an isomorphism.\\

This means that $u$ is injective, $v$ is surjective, so
$$
\deg(V'_\alpha) - \deg(V_\alpha) = \deg (\tau_\alpha) \ge
\deg(\tau)=\deg(V')-\deg(V),
$$
which is exactly the inequality we wanted to prove.\\
 
\begin{eqnarray}\label{ff}\xymatrix{
0 \ar[r] & V \ar[r] & V' \ar[r] & \tau \ar[r] & 0 \\ 0 \ar[r] &
V_\alpha \ar[r] \ar@{^{(}->}[u] & V'_\alpha \ar[r] \ar@{^{(}->}[u] &
\tau_\alpha \ar[r] \ar[u]^{v} & 0 \\  0 \ar[r] & V(-p) \ar[r]
\ar@{^{(}->}[u] & V'(-p) \ar[r] \ar@{^{(}->}[u] & \tau \ar[r]
\ar[u]^{u} & 0 \\  }\end{eqnarray}

\cvd \\

\begin{remark}\label{ennei} \rm Use the notations for the 
multiplicities introduced at the end of the second section, 
to construct the scheme of parabolic quotients. Then the same 
proof shows the following inequalities: 
\begin{eqnarray}\label{var-cellI}
\deg(V_*) \le \deg(V) + \deg(V'_*) - \deg(V') \le \deg(V'_*).
\end{eqnarray}\\
\end{remark}

\bigskip

\subsection{Parabolic extensions} \hfill

\bigskip

Let $\cf '_*, \ \cf ''_*$ be families of parabolic bundles,
parametrized by a scheme $S$. We want to describe a parameter
space for isomorphism classes of nonsplitting parabolic extensions of
type
$$\xymatrix{ 0 \ar[r] & {\cf'_*}_s \ar[r] & F_* \ar[r] & {\cf''_*}_s
\ar[r] & 0,\\ }$$ for $s \in S$. This is actually a consequence of
Lange's results \cite{lange}, so we just introduce the argument
needed to adapt them to the parabolic case.\\

Let $\pi_S$ be the projection $X\times S \epi {} S$ and $R^i{\pi_S}_*
\fhom(\cf ''_*,\cf '_*)$ be the higher direct image sheaves, for
$i=0,1$. For $s \in S$, denote by
$$
\tau ^i_s : R^i{\pi_S}_* \fhom(\cf ''_*,\cf '_*)\ts k(s) \ra
H^i(\fhom({\cf ''_*}_s,{\cf '_*}_s))
$$
the natural base change morphism. The condition that $\tau ^i_s$ is an
isomorphism for all points $s \in S$ will be shortened in $R^i{\pi_S}_*
\fhom(\cf ''_*,\cf '_*)$ {\it commutes with base change}. \\

Let $\mu : S' \ra S$ be a morphism of schemes and denote by
$$
{\rm E}_*(S')=H^0(S',R^1{\pi_{S'}}_* \fhom(\mu ^*\cf ''_*,\mu ^*\cf
'_*)).
$$
Then ${\rm E}_*$ actually is a functor from the category of
$S$-schemes to the category of sets. In fact, let $\nu :S'' \ra S'$ be
a morphism of schemes over $S$. This gives a map ${\rm E}_*(S') \ra
{\rm  E}_*(S'')$ by composition of the natural map
$$
H^0(S',R^1{\pi_{S'}}_* \fhom(\mu ^*\cf ''_*,\mu ^*\cf '_*))  \ra
H^0(S'',\nu ^*R^1{\pi_{S'}}_* \fhom(\mu ^*\cf ''_*,\mu ^*\cf '_*))
$$ 
and the morphism induced in cohomology by the base change morphism
$$
\nu ^*R^1{\pi_{S'}}_* \fhom(\mu ^*\cf ''_*,\mu ^*\cf '_*) \ra
R^1{\pi_{S''}}_* \nu ^*\fhom(\mu ^*\cf ''_*,\mu ^*\cf '_*).
$$
Since it is $\nu ^*\fhom(\mu ^*\cf ''_*,\mu ^*\cf '_*)\cong
\fhom(\nu ^*\mu ^*\cf ''_*,\nu ^*\mu ^*\cf '_*)$, this gives the
morphism ${\rm E}_*(S') \ra {\rm E}_*(S'')$. \\

\begin{proposition} (\cite{lange}, proposition 3.1) Suppose that
  $R^i{\pi_S}_* \fhom(\cf ''_*,\cf'_*)$ commutes with base change,
  for $i=0,1$. Then the functor ${\rm E}_*$ is representable by the
  bundle associated with the locally free sheaf $R^1{\pi_S}_* \fhom(\cf
  ''_*,\cf '_*)^\vee$.\\
\end{proposition}

Let ${\rm PE}_*(S')$ denote the set of invertible quotients
$$
R^1{\pi_S}_*\fhom(\mu ^*\cf ''_*,\mu ^*\cf '_*)^\vee \epi{} \cl.
$$
This defines a functor from the category of
$S$-schemes to the category of sets. \\

\begin{proposition}\label{parext} (\cite{lange}, proposition 4.2)
  Suppose that $R^i{\pi_S}_* \fhom(\cf ''_*,\cf'_*)$ commutes with
  base change, for $i=0,1$. Then the functor ${\rm PE}_*$ is
  representable by the projective bundle $P(R^1{\pi_S}_* \fhom(\cf
  ''_*,\cf '_*)^\vee)$.\\
\end{proposition}

This result is applied in the proof of theorem \ref{gg} in the following
way. Suppose that, for all $s \in S$, there is an isomorphism induced by base
change   
$$
R^1{\pi_S}_* \fhom(\cf ''_*,\cf '_*)\ts k(s)\cong
H^1(\fhom({\cf ''_*}_s,{\cf '_*}_s))
$$ 
and $H^0(\fhom({\cf''_*}_s,{\cf '_*}_s))=0$. 
Then for $i=0,1$ the sheaves $R^i{\pi_S}_*
\fhom(\cf ''_*,\cf '_*)$ commute with base change, the sheaf
$R^1{\pi_S}_* \fhom(\cf ''_*,\cf '_*)$ is locally free over $S$ and
its fibre over a point $s$ is isomorphic to $H^1(\fhom({\cf ''_*}_s,{\cf
'_*}_s))$.  By proposition \ref{parext} the projective bundle
associated with the sheaf $R^1{\pi_S}_* \fhom(\cf ''_*,\cf '_*)^\vee$
parametrizes isomorphism classes of nonsplitting parabolic extensions
of parabolic bundles of the family $\cf ''_*$ by parabolic bundles of
the family $\cf '_*$.\\

\bigskip

\subsection{Proof of theorem \ref{gg}} \hfill

\bigskip

We first prove the theorem for generic weights of $\partial_0 W$:
suppose that the weights $((d_j(p))_{p \in I},k)$ do not lie on any
Seshadri wall.\\

Consider the stratification of $\ce _{r''}$ given by the
quasi-parabolic invariants of the images of parabolic
morphisms. Let $\ce _{((n'_i(p))_{p \in I}),d''}$ be the family of
isomorphism classes of stable parabolic bundles $F_*$, such that there
exists a morphism $\varphi_* : E_* \ra F_*$ for which the vector
bundle $\im(\varphi)=V$ has degree $d''$ and the induced parabolic
structure on $V'$, the saturation of $V$ in $F$, has multiplicities
$((n'_i(p))_{p \in I})$. The parabolic morphism $\varphi_*$ gives rise
to a commutative diagram \\
$$\xymatrix{ 
& & 0 \ar[d] & 0 \ar[d] & \\ 
0 \ar[r] & V \ar[r]\ar@{=}[d] & V'_* \ar[r]\ar[d] & \tau \ar[r]\ar[d]& 0\\ 
0 \ar[r]& V \ar[r] & F_* \ar[r] \ar[d] & G' \ar[r]\ar[d] & 0\\ 
& & G_*\ar@{=}[r]\ar[d] & G \ar[d] & \\ 
& & 0  & 0 & \\ 
}$$ 
and so the exact
sequence in the second column is a nonsplitting parabolic
extension. The bundle $G_*$ has rank $\ell k - r''$, parabolic multiplicities
$$
((\ell d_1(p)-n'_1(p),\dots, \ell d_{l_p}(p)-n'_{l_p}(p))_{p \in I})
$$
and if $t=\deg (\tau)$, then $\deg(G)=\deg(F) -(d''+t)$.\\

\bigskip

Let $\cv'_{((n'_i(p))_{p \in I}),d'',t}$ be the family of isomorphism
classes of parabolic bundles $V'_*$ of degree $d''+t$, multiplicities
$((n'_i(p))_{p \in I})$, such that there exists a stable bundle $F_*$
and a morphism $\varphi_* : E_* \ra F_*$ for which $\im(\varphi)$
generates $V'$ as a subbundle of $F$. Any such bundle is an extension of 
a torsion sheaf $\tau$ of degree $t$ by a quotient bundle $V$ of $E$
and by inequalities \ref{var-cellI} of 
remark \ref{ennei}, if the quotient morphism $\varphi$ induces 
the parabolic structure $V_*$, then it is 
$$
\deg(V_*) \le d'' +\frac{1}{k}\sum_{p\in I}\sum_{i>j}n'_i(p)d_j(p).
$$
Denote by $d''_*$ the right hand side of this inequality. This condition 
implies that the quotient $E \epi{} \im(\varphi)$ is a point of 
a finite union of 
parabolic strata of the scheme $\quot_{r'',d''}(E)$, that we denote 
by $\QI$. More explicitly, it is the union of those strata that correspond 
to functions $s_*$, for which 
$\int_{-1}^0 s_\alpha d\alpha \le d''_*$. By the computation of theorem 
\ref{cong}, its dimension is bounded by
\begin{eqnarray}\label{dimQI}
\dim(\QI) \le r''(r-r'')+r(d''_*-d_{r'',*}),
\end{eqnarray}
where $d_{r'',*}$ is the minimal parabolic degree of 
a rank $r''$ quotient bundle of $E$.\\

Remark that in the generic case it is $\tau\cap I=\emptyset$ and then 
it is enough to consider those quotient morphisms of the stratum 
$\pquot(E_*)$, corresponding to the fixed multiplicities 
$((n'_i(p))_{p\in I})$. 
In any case, we draw the following estimate for the dimension of the family:
$$
\dim (\cv'_{((n'_i(p))_{p \in I}),d'',t}) \le 
r''(r-r'')+r(d''_*-d_{r'',*})+ tr''.
$$\\

\bigskip

Let $\cg_{((n'_i(p))_{p \in I}),d'',t}$ be the family of isomorphism
classes of parabolic bundles $G_*$ of rank $\ell k -r''$, degree $\ell
k(\mu(E_*)+g-1)-(d''+t)$ and multiplicities $((\ell
d_i(p)-n'_i(p))_{p \in I})$, which are parabolic quotients of a stable
bundle $F_*$ by a bundle $V_* \in \cv'_{((n'_i(p))_{p \in
I}),d'',t}$. Consider the family of isomorphism classes of the underlying 
vector bundles and denote it by $\cg$. This family is bounded. 
In fact, since any bundle $G$ of $\cg$ is quotient of some parabolic stable 
bundle $F_*$, if we consider a rank $n$ quotient bundle $G \epi{} H$, there 
is a constant $h(n)$ such that $\mu (H) \ge h(n)$. This condition 
ensures the boundedness of $\cg$. Thus there exists a scheme $S$ 
and a vector bundle $\ch$ on $X \times S$ such that, for all $G$ of the 
family $\cg$ there is an isomorphism $G \cong \ch _s$, for some $s \in S$.
By \cite{newstead}, lemma
4.1, we can suppose that the generic bundle of the family $\cg$ is
semistable, \ie $\dim(\cg) \le (\ell k-r'')^2 (g-1) + 1$.\\
For all $p \in I$, let $\iota_p :\{p\} \times S \mono {}X \times S$ denote 
the inclusion morphism and consider the bundle in flag varieties
$$
\cf_p=\Flag_{(\ell d_1(p)-n'_i(p), \dots, \ell
d_{l_p}(p)-n'_{l_p}(p))}(\iota_p^*\ch).
$$
Let $\cf$ denote the fibred product over $X \times S$ of the 
bundles $\cf_p$. Recall that its dimension is given by
$$
\dim(\cf) = \dim(S) + \sum_{p\in I}\sum_{i>j}(\ell d_i(p)-n'_i(p)) 
(\ell d_j(p)-n'_j(p)).
$$
This family parametrizes quasi-parabolic bundles, whose underlying vector
bundle is isomorphic to $\ch_s$, for some $s \in S$. Thus 
$\cg_{((n'_i(p))_{p \in I}),d'',t}$ is a bounded family and moreover
it is
$$
\dim(\cg_{((n'_i(p))_{p \in I}),d'',t})\le (\ell k-r'')^2 (g-1) +1+
\sum_{p\in I}\sum_{i>j}(\ell d_i(p)-n'_i(p)) (\ell d_j(p)-n'_j(p)).\\
$$

Let $\ce _{((n'_i(p))_{p \in I}),d'', t}$ denote the family of isomorphism
classes of stable parabolic bundles $F_*$, which are parabolic
extensions of a bundle  $G_* \in \cg_{((n'_i(p))_{p \in I}),d'',t}$ by
a bundle $V'_* \in \cv'_{((n'_i(p))_{p \in I}),d'',t}$. \\

\bigskip

\begin{lemma} Let $F_*$ be a stable parabolic bundle and 
$$\xymatrix{ 0 \ar[r] & F'_* \ar[r]^{i_*} &  F_* \ar[r]^{p_*} & F''_*
\ar[r] & 0 \\  }$$ a parabolic extension. Then $\Hom(F''_*,F'_*)=0$.\\
\end{lemma}

\pf If there were a nonzero parabolic morphism $\varphi_* : F''_* \ra
F'_*$, there would be an endomorphism of $F_*$, that is 
$i_*\varphi_* p_* : F_* \ra F_*$, which is not a multiple of the 
identity. \cvd\\

\bigskip

From this lemma it follows that $h^0(\fhom(G_*, V'_*))=0$ and so the
dimension of $H^1(\fhom(G_*, V'_*))\cong \ext(G_*,V'_*)$ is constant
for all $V'_*$ of $\cv'_{((n'_i(p))_{p \in I}),d'',t}$ and $G_*$ of
$\cg_{((n'_i(p))_{p \in I}),d'',t}$. Therefore we can compute the
dimension of the first cohomology group as the opposite of the Euler
characteristic:
\begin{eqnarray*}
\begin{split}
\lefteqn{ \dim(\ext(G_*,V'_*))} \\  & = \sum_{p\in I} \sum_{i>j}(\ell
d_i(p)-n'_i(p))n'_j(p) -(\ell k -r'')(d''+t)\\  & \quad +r'' (\ell k
(g-1 + \mu (E_*)) -(d''+t)) +r''(\ell k -r'')(g-1)  \\ & = - \ell k
(d'' +t) +2r''\ell k (g-1) -{r''}^2 (g-1)+r'' \ell k \mu(E_*)\\ &
\quad  + \sum_{p\in I} \sum_{i>j}(\ell d_i(p)-n'_i(p))n'_j(p).
\end{split}
\end{eqnarray*}\\
Proposition \ref{parext} then gives the following bound for the
dimension of the family of extensions 
\begin{eqnarray*}
\lefteqn{ \dim(\ce _{((n'_i(p))_{p\in I}),d'',t})}\\  & \le 
\dim(\cv'_{((n'_i(p))_{p \in I}),d'',t})+\dim(\cg_{((n'_i(p))_{p \in
I}),d'',t})+h^1(\fhom(G_*,V'_*)) -1.
\end{eqnarray*}
The computation then goes as follows:
\begin{eqnarray*}
\begin{split}
\lefteqn{ \dim(\ce _{((n'_i(p))_{p\in I}),d'',t})}\\  
& \le \dim(\QI)+
r''t + (\ell k -r'')^2 (g-1) \\ & \quad + \sum_{p\in I}\sum_{i>j}(\ell
d_i(p)-n'_i(p))(\ell d_j(p)-n'_j(p)) +1 - \ell k (d''+t) +
r'' \ell k \mu (E_*) \\ 
& \quad +
r'' (2\ell k -{r''})(g-1) + \sum_{p\in I} \sum_{i>j}(\ell
d_i(p)-n'_i(p))n'_j(p) -1\\ 
& =  (\ell k)^2(g-1) +\dim(\QI) +
\sum_{p\in I} \sum_{i>j}\ell d_i(p) \ell d_j(p)\\  
& \quad -
\sum_{p\in I} \sum_{i>j}n'_i(p)\ell d_j(p) - \ell k d'' + t(r''-\ell
k) + r'' \ell k \mu (E_*)\\ 
& =  (\ell k)^2(g-1) + \sum_{p\in I}
d_{\ell d_1(p), \dots , \ell d_{l_p}(p)} + t(r''-\ell k)+ \dim(\QI)\\
& \quad - \ell (k d'' +\sum_{p\in I}
\sum_{i>j}n'_i(p)d_j(p)-\frac{r''}{r}\sum_{p\in I}
\sum_{i>j}n_i(p)d_j(p)).
\end{split}
\end{eqnarray*}\\

The right hand side of the inequality should be read as
\begin{eqnarray*}
\begin{split}
\lefteqn{\dim(M_\ell^{'par}) -1 +t(r''-\ell k)+ \dim(\QI)} \\ 
& \quad - \ell (k d'' +\sum_{p\in I}
\sum_{i>j}n'_i(p)d_j(p)-\frac{r''}{r}\sum_{p\in I}
\sum_{i>j}n_i(p)d_j(p)).
\end{split}
\end{eqnarray*}
By assumption it is $\ell k \ge r \ge r''$ and
$t\ge 0$, so to prove the theorem it is enough to show that
\begin{eqnarray*}
\lefteqn {\dim(\QI)\le} \\  
& \displaystyle \ell (kd''+ \sum_{p\in
I}\sum_{i>j}n'_i(p)d_j(p) - \frac{r''}{r} \sum_{p\in
I}\sum_{i>j}n_i(p) d_j(p)). \notag
\end{eqnarray*}
We can rewrite the right hand side as $\ell \, k(d''_*-r''\mu(E_*))$. Thus,
by inequality \ref{dimQI} it will be enough to show that
$$
r''(r-r'') + r(d''_*-d_{r'',*}) \le \ell \, k(d''_*-r''\mu(E_*)), 
$$
that we can rewrite as
\begin{eqnarray}\label{disug}
r''(r-r'')\le (\ell \, k-r)(d''_*-d_{r'',*}) +\ell \, k(d_{r'',*}
-r''\mu(E_*)).
\end{eqnarray}\\

\bigskip

By assumption it is $\ell \ge \frac{r}{k}$ and by remark 
\ref{ennei} $d''_*$
is greater than or equal to the minimal parabolic degree $d_{r'',*}$.
So in order to get inequality \ref{disug} it is enough to show that 
$$
r''(r-r'') \le \ell \, k(d_{r'',*}-r''\mu(E_*)).
$$\\

This inequality is trivial, when $r''=r$, since in this case both sides 
are equal to zero.\\

\begin{remark} \label{ell} \rm
Suppose $r''$ is strictly less than $r$. Then it is $\ell \le 
\ell \, k(d_{r'',*}-r''\mu(E_*))$. In fact, let $V$ be a 
quotient bundle of rank $r''$ and minimal parabolic degree $d_{r'',*}$.  
The level $k$ of the parabolic structure is such
that $k \mu (E_*) \in \Z$ and moreover $kr''\mu(V_*) = k \deg(V_*)$ 
is an integer as well. Hence the difference $k r'' (\mu (V_*)-\mu
(E_*))$ is an integer, which is strictly positive since $E_*$ is stable. 
Then we draw the inequality
$$
\ell \le \ell \ kr''(\mu(V_*)-\mu(E_*)) = \ell \, k(d_{r'',*}-r''\mu(E_*)).\\
$$
\end{remark}

\bigskip

By this remark inequality \ref{disug} for the nontrivial case
$r'' < r$ follows from the assumption $\ell \ge r''(r-r'')$.
This finishes the proof of theorem \ref{gg} for generic weights.\\

\bigskip

We are left with the case in which the weights $((d_j(p))_{p \in
  I},k)$ of the parabolic structure are on a Seshadri wall and the
  bundle $E_*$ is strictly semistable. Let $E_*={E_0}_* \supset
  {E_1}_* \supset \dots \supset {E_n}_* \supset 0$ denote a
  Jordan-H{\"o}lder filtration of $E_*$, that is each quotient
  ${E_h}_*/{E_{h+1}}_*$ is a stable bundle of parabolic slope
  $\mu(E_*)$. By the previous computation, for each $h$ there is an
  open subscheme $\cu _h$ of the moduli space $M_\ell^{'par}$ 
such that, for all
  stable bundle $F_*$ whose isomorphism class is in $\cu _h$, it is
  $\Hom({E_h}_*/{E_{h+1}}_*, F_*)=0$. Since the moduli space $M_\ell^{'par}$
is irreducible, the open subscheme
$\cu = \cap_{h=1}^l \cu _h$ is nonempty and by definition,
for all stable bundle
  $F_*$ whose isomorphism class is in 
$\cu$ it is $\Hom(E_*, F_*)=0$.\\  
This finishes the proof of the theorem.\\

\bigskip

This bound for the order of base point freeness does 
not depend on the degree $|I|$ of the parabolic divisor and extends 
the result of Popa and Roth for the classical case as well as the result 
of Pauly for rank $2$ parabolic bundles with generic parabolic divisor
of small degree. \\


\bigskip

\begin{thebibliography}{2}

\bigskip

\bibitem[B-H]{boden-hu} Boden, Hans U. and Hu, Yi, {Variations of
moduli of parabolic bundles}, {\it Math. Ann.}, 301: 539--559, 1995.

\bibitem[B-Y 1]{boden-yokogawa1} Boden, Hans U. and Yokogawa,
K{\^o}ji, {Moduli spaces of parabolic {H}iggs bundles and parabolic
${K}({D})$ pairs over smooth curves. {I}}, {\it Internat. J. Math.}, 7: 573--598, 1996.

\bibitem[B-Y 2]{boden-yokogawa2} Boden, Hans U. and Yokogawa,
K{\^o}ji, {Rationality of moduli spaces of parabolic bundles}, {\it
J. London Math. Soc. (2)}, 59: 461--478, 1999.

\bibitem[BP-Gr-Ne]{newstead} Brambila-Paz, L. and Grzegorczyk, I. and
Newstead, P. E., {Geography of {B}rill-{N}oether loci for small
slopes}, {\it J. Algebraic Geom.}, 6: 645--669, 1997.

\bibitem[F]{faltings} Faltings, Gerd, {Stable ${G}$-bundles and
projective connections},{\it J. Algebraic Geom.}, 2: 507--568, 1993.

\bibitem[Ln]{lange} Lange, Herbert, {Universal families of
extensions}, {\it J. Algebra}, 1: 101--112, 1983.

\bibitem[La-So]{laszlo-sorger} Laszlo, Yves and Sorger, Christoph, {The
line bundles on the moduli of parabolic ${G}$-bundles over curves and
their sections}, {\it Ann. Sci. {\'E}cole Norm. Sup. (4)}, 30:
499--525, 1997.

\bibitem[LP]{jlp} Le Potier, Joseph, {Module des fibr{\'e}s
semi-stables et fonctions th{\^e}ta}, Lecture Notes in Pure and
Appl. Math., 179, {\it Dekker}: 83--101, 1996.

\bibitem[Ma]{manivel} Manivel, Laurent, {Fonctions sym{\'e}triques,
polyn{\^o}mes de {S}chubert et lieux de d{\'e}g{\'e}n{\'e}rescence},
Cours Sp{\'e}cialis{\'e}s 3, {\it Soci{\'e}t{\'e} Math{\'e}matique de
France}, 1998.

\bibitem[M-Y]{maruyama-yokogawa} Maruyama, Masaki and Yokogawa,
K{\^o}ji, {Moduli of parabolic stable sheaves}, {\it Math. Ann.}, 293: 77--99, 1992.

\bibitem[Me-S]{seshadri-mehta} Mehta, V.B. and Seshadri C.S., {Moduli
of vector bundles on curves with parabolic structures}, {\it
Math. Ann.}, 248: 205--239, 1980.

\bibitem[Mu-Sa]{mukai-sakai} Mukai, Shigeru and Sakai, Fumio, {Maximal
subbundles of vector bundles on a curve}, {\it Manuscripta Math.}, 52: 251--256, 1985.

\bibitem[N-R]{narasimhan-ramadas} Narasimhan, M. S. and Ramadas,
T. R., {Factorisation of generalised theta functions. {I}}, {\it
Invent. Math.}, 114: 565--623, 1993.

\bibitem[Pa1]{pauly 1} Pauly, Christian, {Espaces de modules de
fibr{\'e}s paraboliques et blocs conformes}, {\it Duke Math. J.}, 84:
217--235, 1996.

\bibitem[Pa2]{pauly 2} Pauly, Christian, {Fibr{\'e}s paraboliques de
rang 2 et fonctions th{\^e}ta g{\'e}n{\'e}ralis{\'e}es}, {\it
Math. Z.}, 228: 31--50, 1998.

\bibitem[Po]{popa} Popa, Mihnea, {Dimension estimates for {H}ilbert
schemes and effective base point freeness on moduli spaces of vector
bundles on curves}, {\it Duke Math. J.}, 107:  469--495, 2001.

\bibitem[Po-Ro]{popa-roth} Popa, Mihnea and Roth, Mike, {Stable maps
and Quot schemes}, {\it Invent. Math.}, to appear.

\bibitem[S]{seshadri} Seshadri C.S., {Fibr{\'e}s vectoriels sur les
courbes alg{\'e}briques}, Ast{\'e}risque, 96. {\it Soci{\'e}t{\'e}
Math{\'e}matique de France}, 1982.

\bibitem[Si]{simpson} Simpson, Carlos T., {Harmonic bundles on
noncompact curves}, {\it J. Amer. Math. Soc.}, 3: 713--770, 1990.

\bibitem[Y]{yokogawa} Yokogawa, K{\^o}ji, {Infinitesimal deformation
of parabolic {H}iggs sheaves}, {\it Internat. J. Math.}, 6: 125--148,
1995.\\


\end{thebibliography}
\end{document}